\documentclass[11pt]{amsart}
\usepackage{amsfonts}
\usepackage{mathrsfs}
\usepackage{bbm}
\usepackage{mathrsfs}
\usepackage{mathrsfs}
\usepackage{amssymb}
\usepackage{hyperref}
\usepackage{amsmath,amssymb, amsthm}
\newcommand{\SL}{\mathrm{SL}(2,\mathbb{R})}

\def\bk{\mathbf{k}}

\def\ti{\tilde}

\let\newpf\proof \let\proof\relax 
\newenvironment{pf}{\newpf[\proofname]}{\qed\endtrivlist}

\newcommand{\ba}{\overline{A}}

\newcommand{\rot}{\mathrm{rot}}

\def\be{\begin{equation}}
\def\ee{\end{equation}}

\def\ba{{\begin{align}}}
\def\ea{{\end{align}}}

\def\bm{\begin{matrix}}
\def\em{\end{matrix}}

\def\SL{{\mathrm{SL}}}

\def\0{{\mathbf 0}}

\newcommand{\la}{\langle }
\newcommand{\ra}{\rangle }

\newtheorem{Theorem}{Theorem}[section]
\newtheorem{Lemma}{Lemma}[section]
\newtheorem{Proposition}{Proposition}[section]
\newtheorem{Corollary}{Corollary}[section]
\newtheorem{Remark}{Remark}[section]

\numberwithin{equation}{section}

\def \bn {\hfill \\ \smallskip\noindent}

\theoremstyle{definition}

\def\proof{\bn {\bf Proof.} }

\def\ssm{\smallsetminus}

\renewcommand{\setminus}{\ssm}

\newcommand{\C}{{\mathbb C}}

\newcommand{\Q}{{\mathbb Q}}
\newcommand{\R}{{\mathbb R}}
\newcommand{\T}{{\mathbb T}}

\newcommand{\Z}{{\mathbb Z}}

\def\B0{{\bold{0}}}

\catcode`\@=12

\def\Empty{}
\newcommand\oplabel[1]{
  \def\OpArg{#1} \ifx \OpArg\Empty {} \else
    \label{#1}
  \fi}

\newcommand{\comm}[1]{}
\newcommand{\comment}[1]{}

\begin{document}

\title[Local Embedding Theorem]{Embedding
 of Analytic Quasi-Periodic Cocycles into  Analytic Quasi-Periodic Linear Systems and its Applications}

\author {Jiangong You}
\address{
Department of Mathematics, Nanjing University, Nanjing 210093, China
} \email{jyou@nju.edu.cn}

\author{Qi Zhou}
\address{
Department of Mathematics, Nanjing University, Nanjing 210093, China
} \email{qizhou628@gmail.com}

\date{\today}

\begin{abstract}
In this paper, we prove that any analytic quasi-periodic cocycle
close to constant is the Poincar\'{e} map of an analytic
quasi-periodic linear system close to constant. With this local
embedding theorem, we get fruitful new results. We show that the
almost reducibility of an analytic quasi-periodic linear system is
equivalent to the almost reducibility of its corresponding
Poincar\'e cocycle.
 By the local embedding theorem and the equivalence, we  transfer  the recent local almost
reducibility results of quasi-periodic linear systems \cite{HoY} to
quasi-periodic cocycles, and the  global reducibility results of
quasi-periodic cocycles \cite{A,AFK} to quasi-periodic linear
systems. Finally, we give a positive answer to a question of
\cite{AFK} and use it to prove Anderson localization results for
long-range quasi-periodic operator with Liouvillean frequency, which
gives a new proof of \cite{AJ05,AJ08,BJ02}. The method developed in
our paper can also be used to prove some nonlinear local embedding
results.
\end{abstract}

\setcounter{tocdepth}{1}

\maketitle

\section{Motivations and main results}

We are concerned with  smooth quasi-periodic linear systems
\begin{eqnarray}\label{qps}
\left\{ \begin{array}{l}\dot{x}=A(\theta)x \\
\dot{\theta}=\omega,
\end{array} \right.
\end{eqnarray}
where $x\in \mathbb{R}^2$, $\theta\in \mathbb{T}^d$, $\omega\in
\mathbb{R}^d$ is rational independent and $A\in C^r (\mathbb{T}^d,
sl(2, \mathbb{R})),$ $r\in N\cup\{\infty,\omega\},$ we denote it
 by $(\omega,A).$ Typical examples are Schr\"odinger systems
where
\begin{eqnarray*}
A(\theta)=V_{E,q}(\theta)=\left( \begin{array}{ccc}
 0 &  1\cr
 q(\theta)-E &  0\end{array} \right)\in sl(2,\mathbb{R}).
\end{eqnarray*}

The time discrete counterparts of the quasi-periodic linear systems
are smooth quasi-periodic $SL(2,\R)$ cocycles:
\begin{eqnarray*}\label{cocycle}
(\mu,\mathcal{A}):&\T^{d-1} \times \R^2 \to \T^{d-1} \times \R^2\\
\nonumber &(\theta,v) \mapsto (\theta+\mu,\mathcal {A}(\theta) \cdot
v),
\end{eqnarray*}
where $\mu \in \T^{d-1}$ with $(1,\mu)$ being rational independent,
$\mathcal{A}\in C^r(\T^{d-1}, SL(2, \mathbb{R}))$. Typical examples
are quasi-periodic Schr\"{o}dinger cocycles where
\begin{eqnarray*}
\mathcal {A}(\theta)=S_E^V(\theta)=\left( \begin{array}{ccc}
 V(\theta)-E &  -1\cr
  1 & 0\end{array} \right)\in SL(2,\mathbb{R}).
\end{eqnarray*}
They are related  to one-dimensional quasi-periodic Schr\"{o}dinger
operators on $l^2(\mathbb{Z})$:
\begin{equation}\label{schro}
(H_{V,\mu,\phi} x)_n= x_{n+1}+x_{n-1} + V( n\mu  + \phi) x_n= E x_n.
\end{equation}
We denote by $C^r_0(\T^{d-1}, SL(2, \mathbb{R}))$ the set of maps
$\mathcal{A}\in C^r(\T^{d-1}, SL(2, \mathbb{R}))$ that are homotopic
to the identity.

 If $d=2,$ it is also interesting to consider the dual
operator of $(\ref{schro}),$ the Long-range quasi-periodic operator:
\begin{equation}\label{long-range}
(L_{V,\alpha, \varphi}\psi)_n=\sum _{k\in\Z}
V_k\psi_{n-k}+2cos2\pi(\varphi+n\alpha)\psi_n,
\end{equation}
where $\alpha \in \R\backslash \Q,$ $V_k$ are fourier coefficients
of $V(\theta)\in C^r(\T,\R).$ These two operators are closely
related by Aubry duality (consult section \ref{pre-aub} for more
details).

 We say that $(\omega,A)$
is $C^r$ reducible, if there exist $B\in C^r(2\T^d,SL(2,
\mathbb{R}))$ and $A_* \in sl(2,\R)$ such that $B$ conjugates
$(\omega,A)$ to $(\omega,A_*)$. It is clear that the concept of
reducibility defined above for Liouvillean frequencies is too
restrictive since in general even a $\R$-valued cocycle is not
reducible. So we need to introduce the weaker concept \textit{almost
reducibility}. A system
 $(\omega,A)$ is  $C^r$ \textit{almost reducible (resp. almost
rotations reducible)} if there exist  sequences of $B_n\in
C^r(2\mathbb{T}^d, SL(2,\mathbb{R}))$, $A_n\in sl(2, \mathbb{R})
(resp. A_n\in C^r (\T^d, so(2, \mathbb{R})))$ and $F_n\in
C^r(\mathbb{T}^d, sl(2,\mathbb{R}))$, such that $B_n$ conjugate
$(\omega,A)$ to $(\omega,A_n+F_n),$ where  $F_n$ is $C^r$ converging
to zero. Another useful concept is \textit{rotations reducibility}.
We say that $(\omega,A)$ is $C^r$ \textit{rotations reducible}, if
there exist $B\in C^r(2\T^d,$ $SL(2, \mathbb{R}))$ and $A_* \in
C^r(\T^d,so(2,\R))$ such that $B$ conjugates $(\omega,A)$ to
$(\omega,A_*)$. These concepts can be defined similarly for cocycles
$\mathcal{A}\in C^r_0(\T^{d-1}, SL(2, \mathbb{R}))$. We say that the
operator $L_{V,\alpha, \varphi},$ (resp.$H_{V,\alpha,\phi}$)
displays Anderson localization, if it has pure point spectrum with
exponentially decaying eigenfunctions.  We remark that reducibility,
almost reducibility and Anderson localization are important issues
in the study of quasi-periodic linear systems and the spectral
theory of Schr\"{o}dinger operators \cite{A,AJ08}.

\subsection{Classical results review}

\subsubsection{Reducibility of  analytic quasi-periodic linear
systems.}  The earliest result of local reducibility was due to
Dinaburg and Sinai~\cite{DS75}, who showed  that $(\omega,
V_{E,q}(\theta))$ is reducible for ``most" sufficiently large $E$,
where $\omega$ is assumed to  satisfy the classical Diophantine
condition:
$$|\langle k, \omega\rangle| \geq
\frac{\gamma^{-1}}{|k|^{\tau}},\quad 0 \neq k \in \Z^d,$$ and
$\gamma, \tau>1$  are fixed positive constants.
 Here $(\gamma,\tau)$ are called the Diophantine constants of $\omega$.  Later, Eliasson~\cite{E92}  obtained the {\it full measure
reducibility} and the {\it almost reducibility} for $(\omega,
V_{E,q}(\theta))$ with Diophantine $\omega$. The proof is based on a
crucial ``resonance-cancellation'' technique which was initially
developed by Moser and P\"oschel~\cite{MP84}.  It should be noted
that these results are {\it perturbative}. Stronger concept is {\it
non-perburbative reducibility}, which  means that the smallness of
the perturbation does not depend on  the Diophantine constants
$(\gamma,\tau).$

However, few results for quasi-periodic linear systems were obtained
since \cite{E92} until \cite{HoY}. Based on a unified approach (KAM
theory and Floquet theory),  Hou and You \cite{HoY} proved that any
two-frequencies quasi-periodic linear system  close to constant is
always {\it almost reducible} and {\it non-perturbative reducible}.
They further proved that it  is analytically {\it  rotations
reducible} if the rotation number is Diophantine w.r.t $\omega$.

\subsubsection{Anderson localization.} Since non-perturbative
reducibility results for quasi-periodic Schr\"odinger cocycles
strongly depends on  Anderson localization for the dual operator, we
review the  Anderson localization results first. For the almost
Mathieu operator $H_{\lambda cos,\alpha,\phi}$, Jitomirskaya
\cite{J} proved that if $\alpha$ is Diophantine, $\lambda>2,$ then
for a.e. $\phi,$ $H_{\lambda cos,\alpha,\phi}$ has Anderson
localization. In the sequel, Avila and Jitomirskaya~\cite{AJ05}
further proved that if $\lambda> 2e^{16\beta/9},\footnote{The
definition of $\beta$ can be found in section \ref{sec:2.1}.}$ then
$H_{\lambda cos,\alpha,\phi}$ has Anderson localization. Similar
results for long-range operator $L_{\lambda V,\alpha, \varphi}$ were
considered by  Bourgain and Jitomirskaya \cite{BJ02}, Avila and
Jitomirskaya~\cite{AJ08}, and we refer to \cite{B2,BG}
 for results on Schr\"odinger operators.

\subsubsection{Reducibility of  analytic quasi-periodic cocycles.}
 Combining Aubry duality~\cite{AA80}
with Anderson localization results \cite{BJ02}, Puig~\cite{Pui06}
obtained a {\it non-perturbative} extension of Eliasson's results.
Avila and Jitomirskaya~\cite{AJ08} further developed a quantitative
Aubry duality  to prove that {\it almost localization} implies {\it
almost reducibility} for the dual model.

Different from the continuous case,  {\it global reducibility} of
quasi-periodic cocycles (reducibility of a cocycle which is not
necessarily close to a constant)  is being developed rapidly.  Based
on Kotani theory \cite{Ko84,Si83} and renormalization scheme
\cite{Kr01}, Avila and Krikorian \cite{AK06} proved that: if
$\alpha$ is recurrent Diophantine, $\mathcal {A}\in
C^{\omega}_0(\T,SL(2,\R))$, then for Lebesgue $a.e. \varphi \in
[0,1],$ $(\alpha,R_{\varphi}\mathcal {A})$\footnote{
$R_{\varphi}:=\left(
\begin{array}{ccc}
cos2\pi\varphi &  -sin2\pi\varphi\\
sin2\pi\varphi &  cos2\pi\varphi
 \end{array}\right)$.} is
either analytic reducible or non-uniformly hyperbolic. With respect
to Liouvillean frequency,  Fayad and Krikorian~\cite{FK} first
obtained that: if $\alpha$ is irrational, $\mathcal {A}\in
C^{\infty}_0(\T,SL(2,\R))$, then for Lebesgue a.e.$ \varphi \in
[0,1],$ $(\alpha,R_{\varphi}\mathcal {A})$ is either $C^{\infty}$
{\it almost rotations reducible} or non-uniformly hyperbolic. Their
main techniques in the proof  were ``algebraic conjugacy trick" and
renormalization scheme~\cite{AK06}. Later, Avila, Fayad and
Krikorian~\cite{AFK} developed this method and obtained a local
positive measure \textit{rotations reducibility} result.
Subsequently, they further proved that for Lebesgue $a.e. \varphi
\in [0,1],$ $(\alpha,R_{\varphi}\mathcal {A})$ is either analytic
{\it rotations reducible} or non-uniformly hyperbolic.

Recently, Avila~\cite{A09} proposed an authentic ``global theory" of
one-frequency $ SL(2,\R)$ cocycles. Cocycles which are not uniformly
hyperbolic are classified into three classes: supercritical,
subcritical and critical. A central issue in the reducibility
 theory is his  Almost Reducible
Conjecture: ``subcritical implies almost reducibility". If $\alpha$
is exponentially Liouvillean, Avila~\cite{A} proved the conjecture
and obtained a corollary: any one-frequency analytic quasi-periodic
cocycle close to constant is {\it almost reducible}.

\subsubsection{Comparison of continuous case and discrete case.} \
In the continuous case, local almost reducibility result is
completely established recently by Hou and You \cite{HoY}, while
there is no result for global reducibility. Since there is no Aubry
duality~\cite{AA80} in the  continuous case, reducibility can not be
obtained by proving localization or almost localization of the dual
system  as in the discrete case \cite{AJ08,Pui06}. However, there is
a uniform way to prove the local almost reducibility results
directly in the continuous case. Unfortunately, the methods
developed in \cite{HoY} can not be applied to quasi-periodic
cocycles directly.
 Global reducibility results are completely  missing for continuous systems since there is no suitable
renormalization scheme as  in the discrete case.  It is still an
interesting question  whether there is corresponding renormalization
scheme for the quasi-periodic linear systems. We note that \cite{Di}
also used renormalization technique to get reducibility results,
however, the method can only be applied to the local situation and
restricted to Brjuno frequency.

In the discrete case, although various global reducibility results
\cite{A,AFK,AK06,FK,Kr01} were obtained,  local almost reducibility
results are not satisfactory. Firstly, there is no unified and
direct approach to deal with the almost reducibility problem for the
quasi-periodic cocycles even in the local regime. The existing
approach which highly depends on the localization results for the
dual model, works only for Schr\"odinger cocycles.
 Secondly, as pointed by Avila, Fayad and Krikorian \cite{AFK}, the study of parabolic behavior was
missing, also there was no generalization of Eliasson's results
\cite{E92}. Readers may refer to   the section $1.1$ of \cite{AFK}
for more discussions.\\

In this paper, we shall establish  a local embedding theorem, which
serves as a  bridge between
 analytic quasi-periodic linear systems and
quasi-periodic cocycles. With this powerful tool, we can deduce
fruitful  new results. For example, one can exchange the almost
reducibility results of quasi-periodic linear systems and
quasi-periodic cocycles for free and then get many missing results
both for the continuous systems and discrete cocycles. Furthermore,
combining Aubry duality, we deduce Anderson localization results for
the long-range operator with Liouvillean frequency. The proof of the
local embedding theorem is interesting in itself and has further
generalizations.

\subsection{Embedding theorem.}

\subsubsection{A local embedding theorem.}  Let $\mathcal {G}=sl(n,\R),
sp(2n,\R), o(n), u(n), so(n)$, and let $G$ be the corresponding Lie
groups. We consider the following $C^r$ quasi-periodic linear
system:
\begin{eqnarray}\label{qls1}
\left\{ \begin{array}{l}\dot{x}=A(\theta)x \\
\dot{\theta}=\omega=(1,\mu),
\end{array} \right.
\end{eqnarray}
where $A \in C^r(\T^d, \mathcal {G})$, $\mu\in \T^{d-1}$ with
$(1,\mu)$ being  rational independent,
$\theta=(\theta_1,\widetilde{\theta})$,
$\widetilde{\theta}=(\theta_2,\cdots,\theta_d)$. Denote by $\Phi^t$
the flow of $(\ref{qls1})$ defined on ${\T}^d\times \R^n$ which is
of the form:
$\Phi^t(\theta_1,\widetilde{\theta},y)=(\theta_1+t,\widetilde{\theta}+t\mu,\Phi^t(\theta_1,\widetilde{\theta})y)$.
We introduce $\mathcal{A}(\cdot)=\Phi^1(0,\cdot)$ which is clearly
defined on ${\T}^{d-1}$ and  called  the corresponding Poincar\'{e}
cocycle defined by $(\ref{qls1}).$ What we are interested in is the
converse, whether we can embed a given $C^r$
 quasi-periodic cocycle into a $C^r$
quasi-periodic flow? Apparently such a cocycle must be homotopy to
the identity.

If $r\neq \omega,$ using the method of suspension flow, Chavaudret
\cite{CC} proved that any $C^r$ smooth $G-$exponential cocycle
\footnote{It means that the cocycle can be written as
$e^{A(\cdot)}.$} can be $C^r$ embedded into the flow.
Rychlik~\cite{Ry} proved that any $C^r$ smooth $SU(2)$ cocycle can
be $C^r$ embedded into a quasi-periodic flow, the proof strongly
depends on the fact that $SU(2)$ is simply connected. Both methods
can't be applied to the analytic case.

 The aim of this paper is to provide
a new method to prove the local embedding of an analytic
quasi-periodic $G-$ valued cocycle into an analytic quasi-periodic
linear system. For a bounded analytic function $F$ defined on $|Im
\theta|<h$, let $\|F\|_h= \sup_{|Im \theta|<h}\|F\|$. We denote by
$C^\omega_{h}(\T^d,*)$ the set of all these $*$-valued functions
($*$ will usually denote $\R$, $\mathcal {G}$). The main theorem is
the following:

\begin{Theorem}\label{localemb}
Let $h>0$, $\mu\in \T^{d-1}$ with  $(1,\mu)$ being  rational
independent,  $A\in \mathcal {G}$, $G\in
C^\omega_{h}(\T^{d-1},\mathcal {G})$. There exist
 $\epsilon=\epsilon(A,h,|\mu|)>0$,
$c=c(A,h,|\mu|)>0$ such that the quasi-periodic cocycle $(\mu,e^A
e^{G(\cdot)})$ can be analytically embedded into a quasi-periodic
linear system provided that $\|G\|_h=\varepsilon<\epsilon$. More
precisely,
 there
exist $\widetilde{A}\in\mathcal {G},$ $F\in
C^\omega_{h/1+|\mu|}(\T^d,\mathcal {G})$ with $\|F\|_{h/1+|\mu|}\leq
c \varepsilon^{1/2}$,  such that $(\mu,e^A e^{G(\cdot)})$ is the
Poincar\'{e} map of
\begin{eqnarray*}
\left\{ \begin{array}{l}\dot{x}=(\widetilde{A}+F(\theta))x \\
\dot{\theta}=(1,\mu).
\end{array} \right.
\end{eqnarray*}
\end{Theorem}

\begin{Remark}
 The selection of $\widetilde{A}$ and precise estimate on $F$ in  Theorems \ref{localemb} will be given explicitly in the proof. We
emphasize
 that $\widetilde{A}=A,$ $\|F\|_{h/1+|\mu|}\leq c \varepsilon$ if $A$  is
 diagonalizable.
 \end{Remark}

\begin{Remark}
The $C^r$ $(r\neq \omega)$ version of  the local embedding theorem
is also true. Different from \cite{CC},  the local structure  is
preserved with precise estimates.
\end{Remark}

\begin{Remark}
The local embedding theorem is also true when  $(1,\mu)$ is rational
dependent. Although the result is trivial by Floquet theory, we will
prove these results in a unified approach.
\end{Remark}

\begin{Remark}
If $d\geq 3$,  an example by Bourgain \cite{B1} shows that
Eliasson's perturbative reducibility result is optimal. Since the
embedding theorem does  not depend on  any Diophantine condition of
the base dynamics, we do embed some non-uniformly hyperbolic
cocycles into quasi-periodic liner systems.
\end{Remark}

\subsubsection{Global embedding results.} We have partial results for
global embedding of analytic  quasi-periodic cocycles into
quasi-periodic linear systems, which are restricted to one-frequency
cocycle. The following  observation is crucial:

\begin{Proposition}\label{emb-conj}
Let  $\mu\in \T^{d-1}$ with  $(1,\mu)$ being rational independent,
$\mathcal {A}\in$ $C^\omega(\T^{d-1}$, $SL(2,\R))$. Suppose that
$(\mu,\mathcal{A})$ is conjugated to $(\mu,\tilde{\mathcal{A}})$,
and $(\mu,\tilde{\mathcal {A}})$ can be embedded into an analytic
quasi-periodic linear system, then $(\mu, \mathcal {A})$ can also be
embedded into an analytic quasi-periodic linear system.
\end{Proposition}

With the help of  Proposition \ref{emb-conj}, we get the embedding
result of uniformly hyperbolic  and almost reducible cocycles:

\begin{Corollary}\label{alemb}
Let $\alpha \in \R\backslash \Q$,  $ \mathcal {A} \in C^\omega_0(\T,
SL(2,\R))$, then we have the following:
\begin{enumerate}
\item Any uniformly hyperbolic cocycle can be analytically embedded
into a quasi-periodic linear system.
\item If $(\alpha, \mathcal{A})$ is almost
reducible, then $(\alpha,\mathcal {A})$ can be analytically embedded
into a quasi-periodic linear system.
\end{enumerate}
\end{Corollary}

As a result of Corollary \ref{alemb}, we obtain that if a global
cocycle is  reduced to the  local regime, then it can be
analytically embedded into a quasi-periodic linear system, if we
recall results of \cite{A,AFK}, then we have:

\begin{Corollary}
Let $\alpha \in \R\backslash \Q$,  $ \mathcal {A} \in C^\omega_0(\T,
SL(2,\R))$, then we have the following:
\begin{enumerate}
\item If $E=\{\varphi \in [0,1]
|L(\alpha,  R_\varphi  \mathcal {A})=0\}$,
 then for almost every $\varphi \in E,$
$(\alpha, R_\varphi \mathcal {A})$ can be analytically embedded into
a quasi-periodic linear system.
\item Assume furthermore $\beta(\alpha)>0$.  If $(\alpha,\mathcal {A})$ is
subcritical, then $(\alpha,\mathcal {A})$ can be analytically
embedded into a quasi-periodic linear system.
\end{enumerate}
\end{Corollary}

\noindent \textbf{Question:} For any $(\alpha,  \mathcal {A}) \in
\R\backslash \Q \times C^\omega_0(\T, SL(2,\R))$, whether it can be
embedded into an analytical quasi-periodical linear system?

\subsubsection{ Nonlinear local embedding results.}  We  point out that
our proof of local embedding theorem is not restricted to linear
case. For example, we can prove the following nonlinear local
embedding result.

\begin{Theorem}\label{qpf}
Let $\rho>0,$ $r>0,$ $s>0,$ and $\mu\in \T^{d-1}$ with $(1,\mu)$
being  rational independent,
 $f\in C^\omega_{r,s}(\T\times\T^{d-1},\R)$.
There exist $c=c(\rho,r,s)>0,$ $\varepsilon=\varepsilon(\rho,r,s)>0$
such that if $\|f\|_{r,s}\leq \varepsilon,$ then the
quasi-periodically forced circle
 diffeomorphism:
\begin{equation}
\left\{
\begin{array}{ll} & \theta \rightarrow \theta + \rho+
f(\theta,\varphi) \\ & \varphi \rightarrow \varphi + \mu
\end{array} \right.
\end{equation}
can be analytically embedded into a quasi-periodically forced circle
flow
\begin{equation}
\left\{
\begin{array}{ll} & \dot{\theta} =\rho+
g(\theta,\varphi) \\ & \dot{\varphi} = \omega=(1,\mu)
\end{array} \right.
\end{equation}
with the estimate $\|g\|_{\frac{r}{1+\rho},\frac{s}{1+|\mu|}}\leq c
\|f\|_{r,s}.$
\end{Theorem}
%
%The second result is about  the perturbation of linear map:
%\begin{Theorem}
%Let $\mu\in \T^{d-1}$, $A\in \mathcal {G}$ being diagonalizable. We
%consider the real analytic perturbation of linear map:
%\begin{eqnarray*}
%F_\varepsilon:& \T^{d-1} \times \R^{n} \to \T^{d-1} \times \R^{n}\\
%\nonumber &(\theta,x) \mapsto ( \theta+\mu, e^A x+ f(\theta,x)),
%\end{eqnarray*}
%when $f$ is small enough, then it can be  analytically embedded into
%the flow
%\begin{eqnarray*}
%\left\{
%\begin{array}{ll} & \dot{x} =A x+
%g(\theta,x) \\ & \dot{\theta} =(1,\mu).
%\end{array} \right.
%\end{eqnarray*}
%\end{Theorem}

\begin{Remark}
Analytic embedding of nearly integrable symplectic maps into
Hamiltonian system is given by Kuksin and P\"oschel \cite{KP}.
\end{Remark}

 Since
the result can be proved with the same method as in our paper, we
will omit the proof. Consult Remark \ref{qpf-lo} for more
discussions. Readers can find more interesting results on
quasi-periodically forced circle
 diffeomorphism in \cite{KWYZ}.

\subsection{Applications of the embedding theorem: from dynamical side.}
\subsubsection{Equivalence.} By applying the local embedding theorem,
we prove that almost reducibility of  an analytic quasi-periodic
system is equivalent to almost reducibility of its corresponding
Poincar\'e cocycle.

\begin{Theorem}\label{dis-con}
 An analytic  quasi-periodic linear system $(\omega,A)$ is almost reducible (resp. rotations
reducible) if and only if its corresponding Poincar\'e cocycle
$(\mu, \mathcal{A})$ is almost reducible (resp. rotations
reducible).
\end{Theorem}

\begin{Remark}
In \cite{Kr}, Krikorian proved that a quasi-periodic linear system
$(\omega,A)$ is reducible
 if and only if its corresponding Poincar\'e cocycle $(\mu, \mathcal{A})$ is reducible.
\end{Remark}

As applications of local embedding theorem and the equivalence of
almost reducibility, we get  many missing results both for the
continuous systems and discrete cocycles.

\subsubsection{Local reducibility of  analytic quasi-periodic
cocycles. }  Let $\alpha \in \R \setminus \Q,$ $\mathcal {A}\in
C^\omega(\T,SL(2,\R)).$ Is there a full Lebesgue measure subset
$\Lambda(\alpha)$, which is explicitly given by  some Diophantine
condition, such that if $\mathcal {A}$ is sufficiently close to
constant and the rotation number $\rot_f(\alpha,\mathcal{A}) \in
\Lambda(\alpha)$, then $(\alpha,\mathcal{A})$ is rotations
reducible? This question was asked  in \cite{AFK}, it  can be seen
as a generalization of Eliasson's result \cite{E92}. We would like
to give an even stronger result to the question:

%The answer has been given in \cite{HoY} for  the continuous systems:

%
%\begin{Theorem}\label{hy}\cite{HoY}  Let
% $\omega=(1,\alpha)$, $\alpha \in \R\backslash \Q$,
%$h>0,$
% $A\in sl(2,\mathbb{R})$ and $F\in C^\omega_{h}(\T^2,sl(2,\R))$. Then there exists  $\epsilon>0$ depending on $A, h$ but not
%on $\alpha$, such that if $\|F\|_h <\epsilon,$ the following results
% hold:
%\begin{enumerate}
%\item The system $(\omega,A+F(\theta))$ is almost
%reducible.
%\item  If the rotation number is Diophantine w.r.t
%$\omega$, then it is analytically rotations reducible.
%\item  Assume furthermore that
%$h>\beta(\alpha)>0$, then it  is analytically reducible.
%\end{enumerate}
%\end{Theorem}
%
%\begin{Remark}
%In \cite{HoY}, it is shown that if $h>3\beta(\alpha),$ then the
%system is reducible. But such conditions are not optimal. Combining
%algebraic conjugacy trick developed in \cite{FK,YZ1}, one can easily
%prove that the optimal condition is $h>\beta(\alpha).$
%\end{Remark}

\begin{Theorem}\label{full-rot}
For every $\alpha \in \R \setminus \Q,$ $h>0,$ $\mathcal {A}\in
C^\omega_h(\T,SL(2,\R))$. If the rotation number
$\rot_f(\alpha,\mathcal {A})$ is Diophantine w.r.t $\omega$,
$$\|\mathcal
{A}-R\|_h<\widetilde{C}min\{h^{2\chi},1\}e^{-\frac{12\pi
h}{1+\alpha}},$$ for some constant matrix $R$, and
$\widetilde{C},\chi$ are numerical constants.
 Then we have the following:
\begin{enumerate}
\item The cocycle $(\alpha,\mathcal {A})$ is analytically rotations reducible.
\item If $\beta(\alpha)=0,$ then $(\alpha,\mathcal {A})$ is analytically  reducible.
\item Furthermore, if $2\pi h>(1+\alpha)\beta>0,$ then $(\alpha,\mathcal {A})$ is
analytically reducible.
\end{enumerate}
\end{Theorem}

\begin{Remark}
 Combining results of \cite{A} and \cite{AFK}, Avila gives an answer
 to the question. Our results will be more natural, and have nice
 spectral applications.
\end{Remark}

In fact, the following stronger result on the  local Almost
Reducibility Conjecture is also a consequence of local embedding
theorem and results of \cite{HoY}.
\begin{Corollary}\label{larc}
Any one-frequency analytic quasi-periodic $SL(2,\R)$ cocycle close
to constant is  almost reducible.
\end{Corollary}

\begin{Remark}
It is a new proof of Avila-Jitomirskaya's theorems \cite{A,AJ08}.
When $\beta(\alpha)=0$, their proof is based on almost localization
results for Long-range operator \cite{AJ08}; when $\beta(\alpha)>0,$
his proof depends on periodic approximation \cite{A}.
\end{Remark}

\subsubsection{Global reducibility of  analytic quasi-periodic systems.}  As
immediate  corollaries of Theorem \ref{dis-con} and the results of
\cite{A,AFK}, we get some global reducibility results for analytic
quasi-periodic linear systems:

\begin{Corollary}\label{cor}
Let $\omega=(1,\alpha)$ with $\alpha \in \R\backslash \Q$, $A\in
C^\omega(\T^2,sl(2,\R))$.  Then we have the following:
\begin{enumerate}
\item For almost every  rotational
number $\rot_f(\omega,A)$, $(\omega,A)$ is either non-uniformly
hyperbolic or (analytically) rotations reducible.
\item Assume further more $\beta(\alpha)>0$. Then $(\omega,A)$ is almost reducible if it is
subcritical.
\end{enumerate}
\end{Corollary}

\begin{Remark}
By Corollary \ref{cor}, we obtain that the Schr\"odinger conjecture
\cite{MMG} (in the essential support of the absolutely continuous
spectrum, the generalized eigenfunctions are almost surely bounded)
is true for continuous quasi-periodic Schr\"odinger operator. This
result  was first verified in a Liouvillean content for the discrete
case \cite{AFK}.
\end{Remark}

\subsubsection{Density of positive Lyapunov exponents.} As an
application of Corollary \ref{cor}, we have the following result:
\begin{Corollary} Let $\omega=(1,\alpha)$ with $\alpha \in \R\backslash
\Q$. Then there is a dense set of  $A\in C^\omega(\T^2,sl(2,\R))$
(in the usual inductive limit topology) such that the linear
quasi-periodic system $(\omega,A)$ has positive Lyapunov exponent.
\end{Corollary}

The result has been proved in the discrete case \cite{A05,FK}. If
$\alpha\in \R$ is recurrent Diophantine, the result was proved in
\cite{Kr}. With the help of Corollary \ref{cor}, the proof can be
carried over for arbitrary irrational $\alpha$ without change.  It
is still open whether uniformly hyperbolic is dense in the category
of analytic quasi-periodic $\SL(2,\R)$-cocycles which are homotopic
to the identity.

\subsection{Applications of the embedding theorem: from spectral side.}

We consider the following long-range quasi-periodic operator:
$$
(L_{\lambda V,\alpha, \varphi}\psi)_n=\lambda \sum _{k\in\Z}
V_k\psi_{n-k}+2cos2\pi(\varphi+n\alpha)\psi_n.
$$
When $\alpha$ is Liouvilean, due to Gordon's lemma \cite{G}, one
often expect that the operator has singular continuous spectrum. Our
result will be if $\beta(\alpha)$ is positive and finite, then for a
suitable range of $\lambda$ Anderson Localization occurs. This
result gives a new proof of previous work \cite{AJ05,AJ08,BJ02}.
What is interesting is that we obtain these Liouvillean results from
the reducibility side, without any localization method.

\begin{Theorem}\label{anderson}
Let $\alpha\in\R\backslash \Q$ be such that $\beta(\alpha)<\infty,$
$2\pi h>(1+\alpha)\beta$, $V \in C^\omega_h(\T,\R)$. Then there
exists a set $\Phi \subset \T$ of full (Lebesgue) measure, such that
if $\phi \in\Phi$, $$\lambda
<\widetilde{C}min\{h^{2\chi},1\}e^{-\frac{12\pi
h}{1+\alpha}}{\|V\|_{(1+\alpha)\beta/2\pi}}^{-1},$$ where
$\widetilde{C},\chi$ are numerical constants,  then the long-range
operator $L_{\lambda V,\alpha,\phi}$ has Anderson Localization.
\end{Theorem}

\begin{Remark}
It is obvious that there exist constant $c_1=c_1(V),$ $c_2=c_2(h),$
such that if $\lambda<e^{-c_1\beta-c_2},$ then $L_{\lambda
V,\alpha,\phi}$ has Anderson Localization.
\end{Remark}

In case $\alpha$ is Diophantine, the result is due to
Bourgain-Jitomirskaya \cite{BJ02}. Avila and Jitomirskaya
\cite{AJ08} proved that if $\alpha$ is not super-liouvillean:
$\xi=\sup_{n>0} \frac{\ln q_{n+1}}{ q_n}<\infty,$ then for
$\lambda<\lambda_0(h,V),$ $L_{\lambda V,\alpha,\phi}$ is almost
localized for all $\phi,$ and has Anderson Localization for
a.e.$\phi$. Actually, though Aubry duality, we can obtain  almost
localization by almost reducibility. We don't pursue this way in
this paper.

If we are restricted to almost Mathieu operator, we can obtain
better estimate:
\begin{Theorem}\label{mathieu}
Let $\alpha\in\R\backslash \Q$ be such that $\beta(\alpha)<\infty.$
If $\lambda<ce^{-\beta},$   with  small constant $c$, then the
almost Mathieu operator $L_{\lambda cos,\alpha,\phi}$
 displays Anderson localization for almost every $\phi$.
\end{Theorem}

\begin{Remark}
In \cite{AJ05}, the authors proved that if
$\lambda<\frac{1}{2}e^{-16\beta/9},$ then if $\beta$ is large, our
result improves theirs. It is still open whether the optimal
condition is $\lambda<\frac{1}{2}e^{-\beta}$ \cite{AJ05}.
\end{Remark}

\noindent \textit{Outline of the paper.}   We first include some
preliminaries  in section $2$. In section $3,$ we give the proof of
the local embedding theorem  and show some global embedding results
in section $4$. As applications, we prove the equivalence of almost
reducibility between the continuous flow and discrete Poincar\'{e}
cocycle in Section $5$. The proof of Theorem \ref{anderson} is shown
in Section $6.$

\section{Preliminaries}

\subsection{Continued Fraction Expansion}\label{sec:2.1}
Let $\alpha \in (0,1)$ be irrational. Define $ a_0=0,
\alpha_{0}=\alpha,$ and inductively for $k\geq 1$,
$$a_k=[\alpha_{k-1}^{-1}],\qquad \alpha_k=\alpha_{k-1}^{-1}-a_k=G(\alpha_{k-1})=\{{1\over \alpha_{k-1}}\},$$
We define
$$p_0=0, \qquad p_1=1$$ $$
q_0=1, \qquad q_1=a_1$$ and inductively,
$$p_k=a_kp_{k-1}+p_{k-2}$$ $$q_k=a_kq_{k-1}+q_{k-2}.$$
It is easy to verify that
$$ \forall 1 \leq k
< q_n,\quad \|k\alpha\|_{\T} \geq \|q_{n-1}\alpha\|_{\T},
$$
and
$$
\|q_n \alpha \|_{\T} \leq {1 \over q_{n+1}}.
$$
thus  $(q_n)$  is the sequence of  denominators of the  best
rational approximations of $\alpha$.

We also denote
$$\beta(\alpha):=\limsup_{n\rightarrow \infty}\frac{\ln
q_{n+1}}{q_n},$$ which means $\beta(\alpha)$ measures how
Liouvillean $\alpha$ is.

\subsection{The rotation number}

 Denote the flow of $(\ref{qps})$ by
$\Phi^t(\theta)$, then we define the rotation number of
$(\ref{qps})$ by $$\rot_f(\omega, A)=\lim_{t\rightarrow
+\infty}\frac{arg(\Phi^t(\theta)x)}{t},$$ where $0\neq x\in
\mathbb{R}^2$, $arg$ denote the angle. It is well-defined and
independent of $(\theta,x)$ \cite{JM82}. The rotation number can be defined
similarly for quasi-periodic cocyles $(\alpha,\mathcal {A})\in
\R\setminus\Q \times C^\omega_0(\T, SL(2, \mathbb{R}))$ \cite{Kr}.
The rotation number $\rot_f$ is said to be rational w.r.t. $\alpha$
if $\rot_f=\frac 12 \langle k_0, \alpha\rangle $ for some $k_0\in
\mathbb{Z}$. It is said to be Diophantine w.r.t. $\alpha$ with some
constants $\gamma,\tau>0$, if $$\|\langle k, \alpha\rangle
-2\rot_f\|_{\R/\Z}\geq \frac{\gamma}{|k|^\tau},\quad 0\neq k\in
\mathbb{Z},$$ and we use  $DC_{\alpha}(\gamma,\tau)$ to denote the
set of all such $\rot_f$. The rotation number is not invariant under
conjugation, but one has the following, the proof can be found in
\cite{Kr}.

\begin{Proposition}\label{rot-conj}
Let $(\alpha,\mathcal {A}_1), (\alpha, \mathcal {A}_2) \in \R
\setminus \Q\times C^r_0(\T, SL(2, \mathbb{R}))$ be two conjugated
quasi-periodic cocycles.  If the conjugacy $B\in C^r(\T, SL(2,
\mathbb{R})$ and has degree $k$, then
$$
\rot_f(\alpha, \mathcal {A}_1)= \rot_f(\alpha, \mathcal {A}_2) +
\langle k, \alpha \rangle \quad \mbox{ mod }1,
$$
If the conjugacy $B\in C^r(2\T, SL(2, \mathbb{R})$ and has degree
$k$, then
$$
\rot_f(\alpha, \mathcal {A}_1)= \rot_f(\alpha, \mathcal {A}_2)
+\frac{1}{2} \langle k, \alpha \rangle \quad \mbox{ mod } 1,
$$
\end{Proposition}

\subsection{Aubry Duality}\label{pre-aub}

Suppose that the eigenvalue equation $H_{V,\alpha,\phi} x=Ex$ has an
analytic quasi-periodic Bloch wave, which means there exist $
\overline{\psi}\in C^\omega(\T, \C)$ and $\varphi \in [0,1)$ such
that
\begin{equation}
x_n = e^{2\pi i n\varphi}  \overline{\psi}\left(n\alpha + \phi
\right).
\end{equation}
We call $\varphi$ the Floquet exponent. If we write $\overline{
\psi}=\sum_{n\in \Z}\psi_n e^{2\pi i n \theta},$ then  direct
computation shows that
\[
\sum_{k \in \Z} V_k \psi_{n-k} + 2 \cos2\pi\left(n\alpha + \varphi
\right) \psi_n = E \psi_n, \qquad n \in \Z,
\]
which means $(L_{V,\alpha, \varphi}\psi)_n=E\psi_n.$  If $\alpha$ is
irrational, then there exists $\sigma^{L}(\lambda
V,\alpha)\subseteq\R,$ such that
$$\sigma^{L}(\lambda V,\alpha)=Spec(L_{\lambda V,\alpha,
\varphi}), \quad \forall \varphi.$$ The rigorous version of the
Aubry duality can be found in \cite{BJ02,GJLS,Pui06}.

Let $\sigma^L_{pp}(V,\alpha,\varphi)$ be the set of point
eigenvalues of $L_{V,\alpha,\varphi}$ which has exponentially
decaying eigenfunctions, and  let $B_{V,\alpha,\varphi}$ be the set
of spectrum of $H_{V,\alpha,\phi}$ which has quasi-periodic Bloch
wave with Floquet exponent $\varphi.$

\begin{Lemma}\label{aubry}
The following facts hold:
\begin{itemize}
\item $\sigma^L_{pp}(V,\alpha,\varphi)=B_{V,\alpha,\varphi}$,
\item $\sigma^{H}(V,\alpha)=\sigma^{L}(V,\alpha).$
\end{itemize}
\end{Lemma}
The proof can be found in \cite{Pui06}.

\section{Proof of Theorem \ref{localemb}}

The local embedding theorem will be proved by Implicit Function
Theorem. The crucial points are the solution of the homological
equation and the construction of suitable Banach spaces. We remark
that our proof doesn't use the method of suspension flow or the
typical property of the Lie group. The method  can also be used to
prove the nonlinear version of the local embedding theorem.

\subsection{Resonance
sites.}

For any $\bk=(k_1,k_2,\cdots,k_d)\in \Z^d$, we define the norm of
$\bk$ by
$$|\bk|=|k_1|+|k_2|+\cdots+|k_d|,$$
and for any $\mu=(\mu_1,\mu_2,\cdots,\mu_d)\in \T^d,$ we define its
norm by
$$|\mu|=|\mu_1|+|\mu_2|+\cdots+|\mu_d|.$$

If $ f(\theta)=\sum_{\bk\in \Z^d}\widehat f(\bk)e^{2\pi  i
\langle\bk,\theta\rangle}\in C_h^\omega(\T^d,\C),$ we use the
weighted norm
$$\|f\|_h:=\sum_{\bk\in \Z^d}|\widehat f(\bk)|e^{2\pi |\bk|h}<\infty,$$
 where
$\theta=(\theta_1,\cdots,\theta_d)$ and
$\langle\bk,\theta\rangle=k_1\theta_1+\cdots +k_d\theta_d.$

 For any $\rho \in
\R$, $\mu \in \T^{d-1}$, $\bk \in \Z^{d-1},$  we  define
$\tilde{k}(\bk) \in \Z$ by
\begin{equation}\label{sel1}|\langle \bk,\mu \rangle+2\rho- \tilde{k}|=\inf_{k\in
\Z}|\langle \bk,\mu \rangle+2\rho-k|.\end{equation} Thus
$\tilde{k}(\bk)$  is uniquely defined if $\inf_{k\in \Z}|\langle
\bk,\mu \rangle+2\rho-k|\ne \frac{1}{2}.$ In case that  $\inf_{k\in
\Z}|\langle \bk,\mu \rangle+2\rho-k|=\frac{1}{2},$
 we choose $\tilde{k}(\bk)$
to be the smaller one which satisfies (\ref{sel1}).

% If $\inf_{k\in
%\Z}|\langle \bk,\mu \rangle+2\rho-k|=0$, that is $\rho$ is rational
%w.r.t. $\omega,$ we have  $\tilde{k}(\bk)=\langle \bk,\mu
%\rangle+2\rho.$ We note it includes the case that  $(1,\mu)$ is
%rational dependent, which means $0$ is rational w.r.t. $\omega.$

 By the construction,  $\tilde{k}(\bk)$  is
uniquely defined and
$$
\tilde{k}(\bk)\in\{[\langle \bk,\mu \rangle+2\rho]-1,[\langle
\bk,\mu \rangle+2\rho],[\langle \bk,\mu \rangle+2\rho]+1\},
$$
where $[\cdot]$ denotes the integer part.

Define the resonance sites $\mathcal {S}_\rho^\mu \subset \Z^{d}$ as
follows
\begin{eqnarray}\label{res set}
\mathcal {S}_\rho^\mu:=\{(\tilde{k}(\bk),\bk), \bk\in \Z^{d-1}\}.
\end{eqnarray}
For any  $f\in  C_h^\omega(\T^{d},\C)$  supported on $\mathcal
{S}_\rho^\mu$, we define its weighted norm by
$$\|f\|_{\rho,h}^\mu:=\sum_{\bk\in \Z^{d-1}}|\widehat
f(\tilde{k}(\bk),\bk)|e^{2\pi |\bk|(1+|\mu|)h},$$ therefore we
define the linear sub-space $\mathcal {B}_{\rho,h}^\mu(\T^{d},\C)$
of
 $C_h^\omega(\T^{d},\C)$:
$$\mathcal {B}_{\rho,h}^\mu(\T^{d},\C)=\{f\in  C_h^\omega(\T^{d},\C) | Supp \widehat f(k_1,\bk)\subset \mathcal {S}_\rho^\mu \}.$$
 The sub-space $\mathcal {B}_{\rho,h}^\mu(\T^{d},\R)$ of  $C_h^\omega(\T^{d},\R)$ is defined similarly.

\begin{Remark}\label{2.1} In  case that
$d=2$, we have
$$e^{-2\pi h(2-2|\rho|)} \|f\|_{\rho,h}^\mu \leq \|f\|_h \leq e^{2\pi h(2|\rho|+1)}
\|f\|_{\rho,h}^\mu,$$ which means that the norms  $\|\cdot
\|_{\rho,h}^\mu$ and $\|\cdot \|_h$ in $\mathcal {B}_{\rho,h}^\mu$
are equivalent. In case that  $d\geq3$, we only have
\begin{equation}\label{equi-norm}\|f\|_h \leq e^{2\pi h(2|\rho|+1)}
\|f\|_{\rho,h}^\mu.\end{equation}
\end{Remark}

  In the following, we will show that  $\mathcal {B}_{\rho,\frac{h}{1+|\mu|}}^\mu(\T^{d},\C)$ is actually isomorphic to
$C_h^\omega(\T^{d-1},\C)$, hence a Banach space. The spaces will be
used to construct the embedded linear system.

\subsection{Embedding operator.}
 For any $f \in
C_h^\omega(\T^{d},\C)$, $\lambda\in\R,$ $\rho\in \R$, we define the
linear operator $$T_{\lambda+i\rho}:
\mathcal{B}_{\rho,\frac{h}{1+|\mu|}}^\mu (\T^{d},\C)\to
C_h^\omega(\T^{d-1},\C)$$ by
$$T_{\lambda+i\rho}
f(\theta)=\int_0^1 f(t,\theta+t\mu)e^{4\pi( \lambda+i\rho) t}dt.$$

If $\lambda\neq0,$ we have
$$T_{\lambda+i \rho} f(\theta)= \sum_{\bk\in \Z^{d-1}}\widehat f(\tilde k(\bk),\bk)\frac{e^{4\pi\lambda +2\pi i
(\tilde k(\bk)+\langle \bk,\mu \rangle+2\rho)}-1}{4\pi\lambda+2\pi
i( \tilde k(\bk)+\langle \bk,\mu \rangle+2\rho)}e^{ 2\pi i
\la\bk,\theta\ra},$$ where  $(\tilde{k},\bk )\in \mathcal
{S}_\rho^\mu,$ consequently, $$\|T_{\lambda+i\rho} f\|_h=
\sum_{\bk\in \Z^{d-1}} |\widehat{T_{\lambda+i\rho} f}(\bk)|e^{2\pi
|\bk|h} \leq \frac{e^{4\pi\lambda}+1}{4\pi\lambda}
\|f\|_{\rho,\frac{h}{1+|\mu|}}^\mu.$$

 If $\lambda=0$ and $\rho$ is not
rational with respect to $\mu$, which means for any $(\tilde{k},\bk
)\in \mathcal {S}_\rho^\mu,$ $\tilde k(\bk)+\langle \bk,\mu
\rangle+2\rho\neq 0,$ then we have
$$T_{i\rho} f(\theta)= \sum_{\bk\in \Z^{d-1}}\widehat f(\tilde k(\bk),\bk)\frac{e^{2\pi i
(\tilde k(\bk)+\langle \bk,\mu \rangle+2\rho)}-1}{2\pi i( \tilde
k(\bk)+\langle \bk,\mu \rangle+2\rho)}e^{ 2\pi i
\la\bk,\theta\ra},$$ consequently,
\begin{equation}\label{em-non}\|T_{i\rho} f\|_h= \sum_{\bk\in
\Z^{d-1}} |\widehat{T_{i\rho}f}(\bk)|e^{2\pi |\bk|h} \leq
\|f\|_{\rho,\frac{h}{1+|\mu|}}^\mu.\end{equation}

If $\lambda=0$ and  $\rho$ is rational with respect to $\mu,$ or in
particular, $(1,\mu)$ is rational dependent, then there exist
$\tilde{\bk}\in \Z^{d-1},$ $\tilde k(\tilde{\bk})\in \Z,$ such that
$\tilde k(\tilde{\bk})+\langle \tilde{\bk},\mu \rangle+2\rho= 0,$ we
have \begin{eqnarray*}T_{i\rho} f(\theta)&=& \widehat f(\tilde
k(\tilde{\bk}),\tilde{\bk})e^{ 2\pi i \la\tilde{\bk},\theta\ra}\\&+&
\sum_{\bk\in \Z^{d-1},\bk\neq \tilde{\bk} }\widehat f(\tilde
k(\bk),\bk)\frac{e^{2\pi i (\tilde k(\bk)+\langle \bk,\mu
\rangle+2\rho)}-1}{2\pi i( \tilde k(\bk)+\langle \bk,\mu
\rangle+2\rho)}e^{ 2\pi i \la\bk,\theta\ra},\end{eqnarray*} then
$(\ref{em-non})$ still holds. Hence in any case,
  $T_{\lambda+i\rho}$
is  a bounded linear operator.

We just point out that when $\rho=0,$
$$T_{\lambda}:\mathcal{B}_{\rho,\frac{h}{1+|\mu|}}^\mu(\T^{d},\R) \rightarrow  C_h^\omega(\T^{d-1},\R) $$
is a bounded linear operator which  maps real functions to real
functions.

We say that $T_{\lambda+i\rho}$ is  an embedding operator if
$T_{\lambda+i\rho}^{-1}$ is a bounded linear operator. In the
following, we prove that $T_{\lambda+i\rho}$ is a linear operator
which does have a bounded inverse.

\begin{Proposition}\label{tech}
For any $\lambda\in\R,$ $\rho \in \R$, $h>0,$ $\mu\in \T^{d-1}$, we
have
$$T_{\lambda+i\rho}^{-1}:C_h^\omega(\T^{d-1},\C)\rightarrow \mathcal{B}_{\rho,\frac{h}{1+|\mu|}}^\mu(\T^{d},\C) $$
is a bounded linear operator.  When $\rho=0$, we have
 $$T_{\lambda}^{-1}:  C_h^\omega(\T^{d-1},\R) \rightarrow  \mathcal{B}_{\rho,\frac{h}{1+|\mu|}}^\mu(\T^{d},\R).$$
\end{Proposition}
\begin{pf} For any $\varphi \in C_h^\omega(\T^{d-1},\C)$, we write
$\varphi(\theta)=\sum_{\bk\in \Z^{d-1}} \hat{\varphi }(\bk) e^{2\pi
i \la\bk, \theta\ra}.$ We first consider the case that $(1,\mu)$ is
rational independent, three cases are
  distinguished.

\noindent \textbf{Case 1} If $\lambda\neq0,$ then we define
\begin{eqnarray*}
\widehat{f}(k_1,\bk)=\left\{
\begin{array}{ccc}\frac{4\pi\lambda+2\pi i( k_1+\langle \bk,\mu
\rangle+2\rho)}{e^{4\pi\lambda+2\pi i
(k_1+\langle \bk,\mu \rangle+2\rho)}-1}\hat{\varphi}(\bk)  & k_1=-\tilde{k}\\
0 & k_1\neq-\tilde{k}
\end{array} \right.
\end{eqnarray*}
where $(\tilde{k},\bk )\in \mathcal {S}_\rho^\mu.$ \\

\noindent \textbf{Case 2} If $\lambda=0$ and $\rho$ is not rational
with respect to $\mu$. In this case, we  define
$\widehat{f}(k_1,\bk)$ by
\begin{eqnarray}\label{sel}
\widehat{f}(k_1,\bk)=\left\{ \begin{array}{ccc}\frac{2\pi i(
k_1+\langle \bk,\mu \rangle+2\rho)}{e^{2\pi i
(k_1+\langle \bk,\mu \rangle+2\rho)}-1}\hat{\varphi}(\bk)  & k_1=-\tilde{k}\\
0 & k_1\neq-\tilde{k}
\end{array} \right.
\end{eqnarray}
where  $(\tilde{k},\bk )\in \mathcal {S}_\rho^\mu.$\\

\noindent \textbf{Case 3} If $\lambda=0$ and $\rho$ is rational with
respect to $\mu$, which means that there exist
$\widetilde{k}_1\in\Z,$ $\widetilde{\bk}_2\in \Z^{d-1},$ such that
$2\rho=-\widetilde{k}_1-\widetilde{\bk}_2\mu$. For
$\bk=\widetilde{\bk}_2,$ we define
\begin{eqnarray*}
\widehat{f}(k_1,\bk)=\left\{ \begin{array}{ccc}\hat{\varphi}(\bk)  & k_1=\widetilde{k}_1\\
0 & k_1\neq\widetilde{k}_1.
\end{array} \right.
\end{eqnarray*}
Otherwise,  for  $\bk\neq \widetilde{\bk}_2$, we define
$\widehat{f}(k_1,\bk)$ by (\ref{sel}).

If $(1,\mu)$ is rational dependent, the construction  is included in
case $3$.

In any cases, by our construction,
$$f(\theta_1,\theta)=\sum_{\bk\in\Z^{d-1} }\widehat{f}(k_1,\bk)e^{2\pi i(
k_1\theta_1+\la\bk, \theta\ra)}$$ is uniquely defined and it
satisfies $T_{\lambda+i\rho}f(0,\theta)=\varphi(\theta).$  Also from
the construction, one sees that $\widehat{f}(k_1,\bk)$ is supported
on the resonance sites $\mathcal {S}_\rho^\mu.$ We now show
$T_{\lambda+i\rho}^{-1}$ is bounded.

If $\lambda\neq0,$ we have
$$\|f\|_{\rho,\frac{h}{1+|\mu|}}^\mu\le\frac{\pi\sqrt{16\lambda^2+1}}{e^{4\pi\lambda}-1} \|\varphi\|_h,$$
which follows by $(k_1,\bk )\in \mathcal {S}_\rho^\mu$ and the
estimate
$$\Big|\frac{4\pi\lambda+2\pi i(
k_1+\langle \bk,\mu\rangle+2\rho)}{e^{4\pi\lambda+2\pi i
(k_1+\langle \bk,\mu\rangle+2\rho)}-1} \Big|\leq
\frac{\pi\sqrt{16\lambda^2+1}}{e^{4\pi\lambda}-1}.
$$

Otherwise, if $\lambda=0,$ by our selection that $(k_1,\bk )\in
\mathcal {S}_\rho^\mu,$ then we have
$$\Big|\frac{2\pi i(
k_1+\langle \bk,\mu \rangle+2\rho)}{e^{2\pi i (k_1+\langle \bk,\mu
\rangle+2\rho)}-1}\Big|<\frac{\pi}{2},$$ consequently, we have
$$\|f\|_{\rho,\frac{h}{1+|\mu|}}^\mu\le \frac{\pi}{2}\|\varphi\|_h.$$
We conclude that $T_{\lambda+i\rho}^{-1}$ is  a bounded linear
 operator, and hence $\mathcal {B}_{\rho,h}^\mu$ is a
Banach space.

  When
$\rho=0$,  $\varphi$ is real analytic, from the formula for
$\widehat{f}(\tilde{k}(\bk),\bk)$, one has
$$\widehat{f}(\tilde{k}(-\bk),-\bk)=\widehat{f}(-\tilde{k}(\bk),-\bk)=\overline{\widehat{f}(\tilde{k}(\bk),\bk)}.$$
This proves that $f$ is real analytic.
\end{pf}

\begin{Remark}
Similar construction was used by Fayad-Katok-Windor in \cite{FKW}.
\end{Remark}

\begin{Remark}
The $C^r$ $(r\neq \omega)$ version of the proposition is also true
(just check the decay of Fourier coeffients), which can be used to
prove the $C^r$ embedding of quasi-periodic cocycles into
quasi-periodic linear systems.
\end{Remark}

\begin{Remark}\label{qpf-lo}
One can also prove that for any $\varphi\in
C^\omega(\T\times\T^{d-1},\R),$ there exists $f\in
C^\omega(\T\times\T^{d},\R),$ such that

$$\int_0^1 f(\theta+\rho t,t,\phi+t\mu)dt=\varphi(\theta,\phi).$$
This fact can be used to prove the nonlinear local embedding of
analytic quasi-periodically forced circle diffeomorphism into
quasi-periodically forced circle flow(c.f. Theorem \ref{qpf}).

\end{Remark}

For any $A\in sl(2,\R),$ let $L:  C_h^\omega(\T^{d},
sl(2,\R))\rightarrow C_h^\omega(\T^{d-1},sl(2,\R))$ be the operator
\begin{equation}\label{operator}L F=\int_0^1 e^{-As}F(s,\theta+s\mu)e^{As}
ds.\end{equation} In the following, we shall prove that there is a
Banach sub-space $\mathcal {B}$ of $C_h^\omega(\T^{d}, sl(2,\R))$,
depending on $A$,  such that
$L:  \mathcal {B}\rightarrow C_h^\omega(\T^{d-1},sl(2,\R))$ is a linear operator with bounded inverse.\\

We recall $sl(2,\R)$ is the set of $2$ by $2$ matrices with real
coefficients of the form $$\left(
\begin{array}{ccc}
 x &  y+z\\
 y-z &  -x
 \end{array}\right)$$
 where $x,y,z\in \R.$ It is
 isomorphic to $su(1,1)$, matrices of the form
$$\left(
\begin{array}{ccc}
 i t &  \nu\\
\bar{ \nu} &  -i t
 \end{array}\right)$$
 with $t\in \R$, $\nu\in \C$. We denote such a matrix by $\{t,\nu\}.$ The isomorphism between $sl(2,\R)$ and
 $su(1,1)$ is given by $B\rightarrow M B M^{-1}$ where
$$M=\left(
\begin{array}{ccc}
 1 &  -i\\
 1 &  i
 \end{array}\right).$$
 Direct calculation shows that
 $$M\left(
\begin{array}{ccc}
 x &  y+z\\
 y-z &  -x
 \end{array}\right)M^{-1}=\left(
\begin{array}{ccc}
 i z &  x-iy\\
x+iy &  -i z
 \end{array}\right).$$
 Denote $H=\left(
\begin{array}{ccc}
 1&  0\\
 0 & -1
 \end{array}\right),$
$J=\left(
\begin{array}{ccc}
 0 &  1\\
 -1 &  0
 \end{array}\right),$ then we have $R:=MJM^{-1}=\left(
\begin{array}{ccc}
 i  &  0\\
0&  -i
 \end{array}\right).$

For any $\rho \in \R$, $h>0,$ $\mu\in \T^{d-1}$, we  define  Banach
spaces
\begin{eqnarray}
\overline{\mathcal {B}}=\left\{\left(\begin{array}{ccc}
i f &  g\\
\bar{g} & -i f
 \end{array}\right)  \Big|
 f\in \mathcal{B}_{0,\frac{h}{1+|\mu|}}^\mu(\T^d, \R), g\in \mathcal{B}_{-\rho,\frac{h}{1+|\mu|}}^\mu(\T^d, \C)\right
 .\},
\end{eqnarray}
and
\begin{eqnarray}\label{space}
\mathcal {B}=M^{-1}\overline{\mathcal {B}}M.
\end{eqnarray}
 We  point out
that $\overline{\mathcal {B}}\subset C^\omega_h(\T^d, su(1,1))$ and
thus $\mathcal {B}\subset C^\omega_h(\T^d, sl(2,\R)),$ since  by
$(\ref{equi-norm})$, we have
$$\|f\|_{\frac{h}{1+|\mu|}}\leq e^{\frac{2\pi h(2|\rho|+1)}{1+|\mu|}}
\|f\|_{\rho,\frac{h}{1+|\mu|}}^\mu.$$

 As a corollary of Proposition
\ref{tech}, we have the following:

\begin{Corollary}\label{inverse}
For any  $\rho \in \R$, $h>0,$ $\mu\in T^{d-1},$  the linear
operator
$$\overline{L}: \overline{\mathcal{B}}\rightarrow C_h^\omega(\T^{d-1},su(1,1))$$ defined
by
\begin{equation}\label{operator1}\overline{L} F=\int_0^1 e^{-2\pi\rho R s}F(s,\theta+s\mu)e^{2\pi \rho Rs}
ds,\end{equation}  is bounded. Moreover, there exists
$C(\rho,h,|\mu|)>0$ such that
$$\overline{L}^{-1}:C_h^\omega(\T^{d-1},su(1,1))\rightarrow \overline{
\mathcal{B}}$$ is  bounded   with  $\|\overline{L}^{-1}\|\leq
C(\rho,h,|\mu|)$.
\end{Corollary}

\begin{pf} For any  $F=\{f_1,f_2\}\in \overline{\mathcal {B}},$ we have
\begin{eqnarray*}
\overline{L}F&=&\int_0^1\left(
\begin{array}{ccc}
i f_1(s,\theta+s\mu) &  f_{2}(s,\theta+s\mu)e^{-4\pi i\rho s}\\
\bar{f_{2}}(s,\theta+s\mu)e^{4\pi i\rho s }& -i f_{1}(s,\theta+s\mu)
 \end{array}\right)ds\\
 &=&\left(
\begin{array}{ccc}
i T_0 f_1 &  T_{-i\rho} f_{2}\\
T_{i\rho}\bar{f_{2}} & -i T_0 f_{1}
 \end{array}\right).
 \end{eqnarray*}
Therefore, $\overline{L}$ is a bounded linear operator.

For any $G =\{g_1,g_2\}\in C_h^\omega(\T^{d-1},su(1,1)),$ by
Proposition $\ref{tech}$ and $(\ref{equi-norm})$, there exists a
unique $f_{1}\in \mathcal{B}_{0,\frac{h}{1+|\mu|}}^\mu(\T^d,\R),$
such that $T_0 f_1=g_1,$ with the estimate $$\|f_{1}\|_h\leq
e^{\frac{2\pi
h}{1+|\mu|}}\|f_{1}\|_{0,\frac{h}{1+|\mu|}}^\mu\leq\frac{\pi}{2}
e^{\frac{2\pi h}{1+|\mu|}}\|g_{1}\|_h.$$ Again by Proposition
$\ref{tech}$ and $(\ref{equi-norm})$,
 there exists a unique $f_{2}\in \mathcal{B}_{-\rho,\frac{h}{1+|\mu|}}^\mu(\T^d,\C),$ such that
$T_{-i\rho} f_{2}=g_{2},$ with the estimate
$$\|f_{2}\|_h\leq e^{\frac{2\pi h(2|\rho|+1)}{1+|\mu|}} \|f_{2}\|_{-\rho,\frac{h}{1+|\mu|}}^\mu\leq \frac{\pi}{2}e^{\frac{2\pi h(2|\rho|+1)}{1+|\mu|}}\|g_{2}\|_h.$$
We remark that $T_{i\rho} \bar{f_{2}}=\bar{g_{2}}.$

 So there exists a unique  $F\in \overline{\mathcal {B}}$  such that
$\overline{L}{F}=G$ with the estimate $$\|F\|_h \leq C(\rho,h,|\mu|)
\|G\|_h,$$ where $C(\rho,h,|\mu|)=\frac{\pi}{2}e^{\frac{2\pi
h(2|\rho|+1)}{1+|\mu|}}.$  It follows that $\overline{L}^{-1}$
exists and is bounded with $\|\overline{L}^{-1}\|\leq
C(\rho,h,|\mu|)$.
\end{pf}

\begin{Corollary}\label{inverse1}
For any  $\rho \in \R$, $h>0,$ $\mu\in \T^{d-1},$ let $A=2\pi\rho
J,$ then the  linear operator
$$L: \mathcal{B}\rightarrow C_h^\omega(\T^{d-1},sl(2,\R))$$ defined
by $(\ref{operator})$ is bounded. Moreover, there exists constant
$C(\rho,h,|\mu|)>0$ such that
$$L^{-1}:C_h^\omega(\T^{d-1},sl(2,\R))\rightarrow
 \mathcal{B}$$
is  bounded   with  $\|L^{-1}\|\leq C(\rho,h,|\mu|)$.
\end{Corollary}

\begin{pf} It is an immediate corollary of Corollary \ref{inverse},
since the Banach spaces $C_h^\omega(\T^{d},sl(2,\R))$ and
$C_h^\omega(\T^{d},su(1,1))$ are isomorphic by $B\rightarrow M B
M^{-1}$.\end{pf}

If $A$ is hyperbolic, the operator defined by $(\ref{operator})$ is
still bounded, which is the following:

\begin{Corollary}\label{inverse2}
For any  $\lambda \in \R$, $h>0,$ $\mu\in\T^{d-1},$ let
$A=2\pi\lambda H,$
\begin{eqnarray*} \widetilde{\mathcal
{B}}=\left\{\left(\begin{array}{ccc}
f_{11} &  f_{12}\\
f_{21} & - f_{11}
 \end{array}\right)  \Big|
 f_{11},f_{12},f_{21}\in \mathcal{B}_{0,\frac{h}{1+|\mu|}}^\mu(\T^d,\R)\right
 \},
\end{eqnarray*}
then the  linear operator
$$L: \widetilde{\mathcal{B}}\rightarrow C_h^\omega(\T^{d-1},sl(2,\R))$$ defined
by $(\ref{operator})$ is bounded. Moreover
$L^{-1}:C_h^\omega(\T^{d-1},sl(2,\R))\rightarrow \widetilde{
\mathcal{B}}$ is  bounded   with  $\|L^{-1}\|\leq
\frac{\pi\sqrt{16\lambda^2+1}}{e^{4\pi\lambda}-1}e^{\frac{2\pi
h}{1+|\mu|}}$.
\end{Corollary}

\begin{pf} We omit the proof since it is similar with  the proof of
Corollary $\ref{inverse}$. \end{pf}

\subsection{Proof of Theorem \ref{localemb}.}

In order to  make the ideas of proof clearly, we will only prove the
theorem in the group $SL(2,\R),$ and just give an outline of the
proof in other Lie groups. The following embedding theorem also
includes the case that $(1,\mu)$ is rational dependent.

\begin{Theorem}\label{localemb-sl}
Suppose that $\mu\in \T^{d-1}$, $h>0,$ $G\in
C^\omega_{h}(\T^{d-1},sl(2,\R))$, $A\in sl(2,\R)$ being constant.
Then there exist $\epsilon=\epsilon(A,h,|\mu|)>0,$
$c=c(A,h,|\mu|)>0,$ $\widetilde{A}\in sl(2,\R)$ and $F\in
C^\omega_{h/1+|\mu|}(\T^d,sl(2,\R))$ such that the cocycle $(\mu,e^A
e^{G(\cdot)})$ is the Poincar\'e map of
\begin{eqnarray}\label{al-ref1}
\left\{ \begin{array}{l}\dot{x}=(\widetilde{A}+F(\theta))x \\
\dot{\theta}=\omega=(1,\mu)
\end{array} \right.
\end{eqnarray}
provided that $\|G\|_{h}= \varepsilon<\epsilon.$ Moreover, we have
the following
\begin{enumerate}
\item If $A$ is in
the real normal forms $\left(
\begin{array}{ccc}
 \lambda &  0\cr
 0 &  -\lambda\end{array} \right)$ or
  $\left(
\begin{array}{ccc}
0  &  \rho\cr
 -\rho &  0\end{array} \right),$ then $\widetilde{A}=A,$ $\|F\|_{\frac{h}{1+|\mu|}}\leq c \varepsilon$.
 \item If
 $A$ is in the real normal form $\left(
\begin{array}{ccc}
 0 &  1\\
 0 &  0
 \end{array}\right),$ then $\widetilde{A}=\left(
\begin{array}{ccc}
 0 &  0\\
 0 &  0
 \end{array}\right)$ and $\|F\|_{\frac{h}{1+|\mu|}}\leq c \varepsilon^{\frac 12}.$
\end{enumerate}

\end{Theorem}

\begin{pf} \textbf{Case 1.} $A$ is elliptic. Without lose of generality,
we assume $A=2\pi\rho J$, and  define $\mathcal {B}$ as in
$(\ref{space})$. Suppose that $ \Phi^t(\theta)$ is the flow of
$(\ref{al-ref1})$,
$$ \Phi^t(\theta)=e^{At}\Big(I +\int_0^t e^{-A s}F(\theta+s\omega)  \Phi^s(\theta)ds \Big),$$
where $I$ denotes the identity matrix. The cocycle $(\mu, e^A
e^{G(\widetilde{\theta})})$ can be embedded into the linear system
$(\ref{al-ref1}),$ which means $\Phi^1(0,\widetilde{\theta})=e^A
e^{G(\widetilde{\theta})}$, i.e.,
\begin{equation}\label{non}e^{A}\Big(I +\int_0^1 e^{-A s}F(s,\widetilde{\theta}+s\mu)  \Phi^s(0,\widetilde{\theta})ds \Big)=e^A e^{G(\widetilde{\theta)}}.\end{equation}

We construct the nonlinear functional
$$\Psi:\mathcal {B}\times C_h^\omega(\T^{d-1},sl(2,\R))\rightarrow C_h^\omega(\T^{d-1},gl(2,\R))$$ by
$$\Psi(F,G)= I +\int_0^1 e^{-A s}F(s,\widetilde{\theta}+s\mu)
\Phi^s(0,\widetilde{\theta})ds - e^{G(\widetilde{\theta})}.$$
Immediate check shows that $\Psi(0,0)=0,$ and

\begin{eqnarray*}D_F\Psi(F,G)(\widetilde{F})&&=\int_0^1 e^{-A
s}\widetilde{F}(s,\widetilde{\theta}+s\omega)\Phi^s(0,\widetilde{\theta})ds\\
&&+\int_0^1 e^{-A s}F(s,\widetilde{\theta}+s\mu) D_F
\Phi^s(0,\widetilde{\theta})\widetilde{F}(s,\widetilde{\theta}+s\mu)ds.\end{eqnarray*}
Consequently, we have
$$D_F\Psi(0,0)(\widetilde{F})=\int_0^1 e^{-A s}\widetilde{F}(s,\widetilde{\theta}+s\mu)e^{A s} ds.$$
 By
Corollary \ref{inverse1}, $D_F\Psi(0,0): \mathcal {B}\rightarrow
C_h^\omega(\T^{d-1},sl(2,\R))$ is a bounded linear operator with
bounded inverse. Moreover, we have the estimate
$\|D_F\Psi(0,0)^{-1}\|\leq C(\rho,h,|\mu|),$ where $C(\rho,h,|\mu|)$
is defined in Corollary \ref{inverse}.

By Implicit Function Theorem, when $$\|G\|_h \leq
\varepsilon<\frac{1}{C(\rho,h,|\mu|)^2},$$ then there exists  $F\in
\mathcal {B}\subset C^{\omega}_{\frac{h}{1+|\mu|}}(\T^d,sl(2,\R))$
with $\|F\|_{\frac{h}{1+|\mu|}}\leq C(\rho,h,|\mu|)\varepsilon,$
such that the nonlinear functional $\Psi(F,G)=0$ has a solution.
That is to say  $(\mu, e^A e^{G(\cdot)})$ is the Poincar\'e map of
$(\ref{al-ref1})$.\\

\noindent\textbf{Case 2} $A$ is hyperbolic. Without lose of
generality, we assume $A=2\pi\lambda H$. In this case, the proof
goes along the same line as in case
 $1.$ We only need to substitute $\mathcal {B}$ by  $\widetilde{\mathcal
 {B}}$, and Corollary \ref{inverse1} by Corollary \ref{inverse2}.\\

\noindent\textbf{Case 3} $A$ is parabolic. Without lose of
generality, we assume $A=\left(
\begin{array}{ccc}
 0 &  1\\
 0 &  0
 \end{array}\right)$. For any $G\in C_h^\omega(\T^{d-1},sl(2,\R))$ with
$\|G\|_h\leq \varepsilon,$ we set $B=\left(
\begin{array}{ccc}
 \varepsilon^{\frac{1}{4}} &  0\\
 0 &  \varepsilon^{-\frac{1}{4}}
 \end{array}\right)$, then we have
 $$B e^A e^{G(\cdot)} B^{-1}=e^{\left(
\begin{array}{ccc}
 0 &   \varepsilon^{\frac{1}{2}}\\
 0 &  0
 \end{array}\right)}e^{\widetilde{G}(\cdot)},$$
 with $\|\widetilde{G}\|_h\leq  \varepsilon^{\frac{1}{2}}.$
 This  means the cocycle $(\mu, e^A e^{G(\cdot)})$ can be seen as a perturbation of
 $(\mu, I)$, then we apply case $1$ to finish the proof.\end{pf}

\begin{Remark}
One can generalize the proof to other Lie groups without difficulty.
Suppose the constant matrix $A$ is diagonalizable, the proof follows
  case $1$ above. Suppose the constant part $A$ has Jordan blocks,
the proof is similar to the case $3$ above.
\end{Remark}

\section{Global embedding results.}

In this section, we first prove that embedding is conjugacy
invariant (Proposition \ref{emb-conj}), and then apply it to prove
some global embedding results.\\

\textbf{Proof of Proposition \ref{emb-conj}.}

 By the assumption, there exist $B\in
C^\omega(2\T^d,SL(2,\R))$ and the flow
$\ti\Phi^t(\theta_1,\widetilde{\theta})\in
C^\omega(2\T^d,SL(2,\R)),$ such that
$$B(\widetilde{\theta}+\mu)\mathcal {A}(\widetilde{\theta})B(\widetilde{\theta})^{-1}= \tilde{\mathcal{A}}(\widetilde{\theta}),$$
and   $$\tilde{\mathcal{A}}(\widetilde{\theta})=\ti
\Phi^1(0,\widetilde{\theta}).$$ For $(x_1,\widetilde{x})\in {\R}^d,$
$$\Phi^{x_1}(0,\widetilde{x}-x_1\mu)=B(\widetilde{x}-x_1\mu+\mu)^{-1}\ti \Phi^{x_1}(0,\widetilde{x}-x_1\mu)B(\widetilde{x}-x_1\mu)$$
 is well-defined, $2$-periodic in $\widetilde{x}$ and analytic.

Once $\Phi^{x_1}(0,\widetilde{x})$ is given, we define
\begin{equation}\label{flow}
\Phi^{t}(x_1,\widetilde{x})=\Phi^{x_1+t}(0,\widetilde{x}-x_1\mu)\big(\Phi^{x_1}(0,\widetilde{x}-x_1\mu)\big)^{-1}.
\end{equation}
 Immediate check shows that $\Phi^{t}(x_1,\widetilde{x})$ is also $2$-periodic in $x_1$ and
analytic, hence $\Phi^{t}\in C^\omega(2\T^d,SL(2,\R))$. Next we show
$(\ref{flow})$ indeed defines a flow. To prove this, let
$t\rightarrow t+s,$ then we have
\begin{eqnarray*}
\Phi^{t+s}(\theta_1,\widetilde{\theta})&=&\Phi^{\theta_1+t+s}(0,\widetilde{\theta}-\theta_1\mu)\big(\Phi^{\theta_1}(0,\widetilde{\theta}-\theta_1\mu)\big)^{-1}\\
&=&\Phi^{t}(\theta_1+s,\widetilde{\theta}+s\mu)\Phi^{\theta_1+s}(0,\widetilde{\theta}-\theta_1\mu)\big(\Phi^{\theta_1}(0,\widetilde{\theta}-\theta_1\mu)\big)^{-1}\\
&=&\Phi^{t}(\theta_1+s,\widetilde{\theta}+s\mu)\Phi^{s}(\theta_1,\widetilde{\theta}).
\end{eqnarray*}
The second equality follows from the substitution
$\theta_1\rightarrow \theta_1+s$, $\widetilde{\theta}\rightarrow
\widetilde{\theta}+s\mu$.

It is clearly that
$$\Phi^{1}(0,\widetilde{\theta})=B(\widetilde{\theta}+\mu)^{-1}\ti
\Phi^{1}(0,\widetilde{\theta})B(\widetilde{\theta})=B(\widetilde{\theta}+\mu)^{-1}
\tilde{\mathcal
{A}}(\widetilde{\theta})B(\widetilde{\theta})=\mathcal
{A}(\widetilde{\theta}),$$ which means $(\mu,\mathcal{A})$ can be
embedded into the flow $\Phi^{t}(\theta_1,\widetilde{\theta})$.\qed

As an immediate consequence of Proposition \ref{tech} and
Proposition \ref{emb-conj}, we have the following embedding result
for uniformly hyperbolic cocycle.

\begin{Corollary}\label{uh-emb}
Let  $\mu\in \T^{d-1},$  $\omega=(1,\mu)$ is rational independent,
$\mathcal {A}\in C^\omega_0(\T^{d-1},SL(2,\R))$. If
$(\mu,\mathcal{A})$ is uniformly hyperbolic, then it can be
analytically embedded into a quasi-periodic linear system.
\end{Corollary}

\begin{pf}If $(\mu,\mathcal{A})$ is uniformly hyperbolic, then there
 exist $\varphi \in C^\omega(\T^{d-1},\R),$ $B\in C^\omega(2\T^{d-1},SL(2,\R))$ such that
 $$B(\cdot+\mu)\mathcal {A}(\cdot)B(\cdot)^{-1}=\widetilde{\mathcal {A}}(\cdot)=\left(
\begin{array}{ccc}
 e^{\varphi(\cdot)} &  0\\
 0 &   e^{-\varphi(\cdot)}
 \end{array}\right).$$
By Proposition \ref{tech}, there exists $f \in C^\omega(\T^d,\R)$
such that $T_0 f=\varphi,$ which means the quasi-periodic cocycle
$(\mu,\widetilde{\mathcal{A}})$ can be embedded into
\begin{eqnarray*}
\left\{ \begin{array}{l}\dot{x}=\left(
\begin{array}{ccc}
 f(\theta) &  0\\
 0 &  -f(\theta)
 \end{array}\right)x \\
\dot{\theta}=\omega=(1,\mu).
\end{array} \right.
\end{eqnarray*}
 Hence the result follows from Proposition
 \ref{emb-conj}.\end{pf}

\begin{Remark}
By Corollary \ref{uh-emb}, we get another proof of Theorem
\ref{localemb-sl} in the case $A$ is hyperbolic, since  uniformly
hyperbolic is an open condition, if $\|G\|_h$ is small, then $(\mu,
e^A e^{G(\cdot)})$ is  uniformly hyperbolic.
\end{Remark}

\begin{Corollary}
Let $\alpha \in \R\backslash \Q$,  $ \mathcal {A} \in C^\omega_0(\T,
SL(2,\R))$.  If $(\alpha, \mathcal{A})$ is almost reducible, then
$(\alpha,\mathcal {A})$ can be analytically embedded into a
quasi-periodic linear system.
\end{Corollary}

\begin{pf} If $(\alpha,\mathcal{A})$ is almost reducible, then there
exist $B_n\in C_{h_n}^\omega(2{\T}, SL(2,{\R}))$, $A_n\in sl(2,\R)$,
$F_n\in C_{h_n}^\omega({\T}, sl(2,{\R}))$ with $\|F_n\|_{h_n}\leq
\varepsilon_n\rightarrow 0$,  such that
$$\mathcal{A}(\cdot)B_n(\cdot)=B_n(\cdot+\alpha)e^{A_n}e^{F_n(\cdot)}.$$
When $n$ is large enough, $\varepsilon_n \leq \varepsilon,$ where
$\varepsilon=\varepsilon(A_n,h_n,\alpha)$ is defined in Theorem
\ref{localemb}. By Theorem \ref{localemb},  the quasi-periodic
cocycle $(\alpha,e^{A_n}e^{F_n(\cdot)})$ can be analytically
embedded into a quasi-periodical linear system. Therefore the result
follows from Proposition \ref{emb-conj}. \end{pf}

\section{Equivalence of almost reducibility results}

We only give the proof in the $SL(2,\R)$ cocycle case, it can be
generalized to other Lie groups without difficulty. In the
following, we let $\mu\in \T^{d-1}$, $\omega=(1,\mu)$ is rational
independent.

\begin{Lemma}\label{almost1}
Let $h>0,$ $A\in C^\omega_h(\T^d,sl(2,\mathbb{R}))$. If $(\omega,A)$
is almost reducible, then the corresponding Poincar\'{e} cocycle
$(\mu, \mathcal{A})$ is almost reducible.
\end{Lemma}

\begin{pf} If $(\omega,A)$ is almost reducible, then  there exist
$B_n\in C_{h_n}^\omega(2{\T}^d, SL(2,{\R})),$ $A_n\in sl(2,\R),$
$F_n\in C_{h_n}^\omega({\T}^d, sl(2,{\R}))$ such that $B_n$
conjugate $(\omega,A)$ to
\begin{eqnarray}\label{al-ref}
\left\{ \begin{array}{l}\dot{x}=(A_n + F_n(\theta))x \\
\dot{\theta}=\omega=(1,\mu)
\end{array} \right.
\end{eqnarray}with $\|F_n\|_{h_n}\leq \varepsilon_n\rightarrow 0$. Denote by $\Phi^t(\theta)$ the flow induced by $(\omega,A).$
Now we fix $n$ which is large enough.\\

\noindent \textbf{Case 1 : $A_n$ is hyperbolic. } In this case,
$(\ref{al-ref})$ is uniformly hyperbolic, and then $(\omega,A)$ is
uniformly hyperbolic. Consequently,
$(\mu,\mathcal{A})$ is  almost reducible since  it is uniformly hyperbolic.\\

\noindent \textbf{Case 2 : $A_n$ is elliptic. } We write
$A_n=2\pi\rho_n J$, suppose that $\ti \Phi^t(\theta)$ is the
corresponding flow of  $(\ref{al-ref})$, then we have
\begin{equation}\label{al-ref3}\ti \Phi^t(\theta)=e^{2\pi t\rho_n J}\Big(I +\int_0^t e^{-2\pi s\rho_n J}F_n(\theta+s\omega) \ti \Phi^s(\theta)ds \Big).\end{equation}
Denote by $G^t(\theta)=e^{-2\pi t\rho_n J}\ti \Phi^t(\theta)$, then
$$
G^t(\theta)=I +\int_0^t e^{-2\pi s \rho_n J}F_n(\theta+\omega
s)e^{2\pi s\rho_n J}G^s(\theta)ds,
$$
let $g(t)=\|G^t(\theta)\|_{h_n}$, then
$$
g(t) \leq 1 + \int_0^t \|F_n\|_{h_n}g(s)ds.
$$
By Gronwall's inequality,  we have $g(t)\leq e^{\varepsilon_n t}.$

In the equation $(\ref{al-ref3})$, let $t=1$, then we have $\ti
\Phi^1(0,\widetilde{\theta})=e^{2\pi \rho_n
J}(I+\widetilde{F}_n(\widetilde{\theta})),$ with the estimate
$\|\widetilde{F}_n\|_{h_n}\leq \int_0^1\varepsilon_n g(t) dt\leq
2\varepsilon_n.$ Since $B_n$ conjugates $(\omega,A)$  to
$(\ref{al-ref}):$
$$\Phi^t(0,\widetilde{\theta})=B_n(t,\widetilde{\theta}+t\mu)\ti \Phi^t(0,\widetilde{\theta})
B_n(0,\widetilde{\theta})^{-1}.$$ Let $t=1$, we have
$$B_n(0,\widetilde{\theta}+\mu)^{-1}\mathcal{A}(\widetilde{\theta})B_n(0,\widetilde{\theta})=e^{2\pi \rho_n J}(I+\widetilde{F}_n(\widetilde{\theta})),$$
which means that the Poincar\'{e} cocycle $(\mu,\mathcal{A})$ is
almost reducible.\\

\noindent \textbf{Case 3 : $A_n$ is parabolic.} Without lose of
generality, we assume $A_n=\left(
\begin{array}{ccc}
 0 &  1\\
 0 &  0
 \end{array}\right)$, let $B=\left(
\begin{array}{ccc}
 \varepsilon_n^{\frac{1}{4}} &  0\\
 0 &  \varepsilon_n^{-\frac{1}{4}}
 \end{array}\right)$, then $x=By$ transformation $(\ref{al-ref})$ to
\begin{eqnarray*}
\left\{ \begin{array}{l}\dot{x}=\bigg(\left(
\begin{array}{ccc}
0 &  \varepsilon_n^{\frac{1}{2}}\\
 0 & 0
 \end{array}\right)+ \overline{F}_n(\theta)\bigg)x \\
\dot{\theta}=\omega=(1,\mu)
\end{array} \right.
\end{eqnarray*}
with $\|\overline{F}_n\|_{h_n}\leq \varepsilon_n^{\frac{1}{2}}$. It
is then  reduced  to  Case $2$.\end{pf}

The converse of Lemma \ref{almost1} is also true:
\begin{Lemma}\label{almost2}
Let $h>0,$ $A\in C^\omega_h(\T^d,sl(2,\mathbb{R}))$, $(\mu,
\mathcal{A})$ is the corresponding Poincar\'{e} cocycle of
$(\omega,A)$. If $(\mu,\mathcal{A})$ is almost reducible,  then
$(\omega,A)$ is almost reducible.
\end{Lemma}

\begin{pf}  If $(\mu,\mathcal{A})$ is almost reducible, then there exist
$B_n\in C_{h_n}^\omega(2{\T}^{d-1}, SL(2,{\R}))$,
$\widetilde{A}_n\in sl(2,\R)$, $\widetilde{F}_n\in
C_{h_n}^\omega({\T}^{d-1}, sl(2,{\R}))$ with
$\|\widetilde{F}_n\|_{h_n}\leq \varepsilon_n$,  such that
\begin{equation}\label{almost-1}\mathcal{A}(\cdot)B_n(\cdot)=B_n(\cdot+\mu)e^{\widetilde{A}_n}e^{\widetilde{F}_n(\cdot)}.\end{equation}
When $n$ is large enough, by Theorem \ref{localemb}, there exists
$\overline{F}_n\in C_{\frac{h_n}{1+|\mu|}}^\omega({\T^d},
sl(2,{\R}))$, such that the quasi-periodic cocycle $(\mu,
e^{\widetilde{A}_n}e^{\widetilde{F}_n(\cdot)})$ can be embedded into
\begin{eqnarray}\label{al-ref2}
\left\{ \begin{array}{l}\dot{x}=(\widetilde{A}_n + \overline{F}_n(\theta))x \\
\dot{\theta}=(1,\mu).
\end{array} \right.
\end{eqnarray}

Suppose that $ \overline{\Phi}^t(\theta)$ is the corresponding flow
of $(\ref{al-ref2})$, we thus  extend $B_n(\cdot)$ to the torus
$2{\T}^d$ in  the following way: for $(x_1,\widetilde{x})\in
{\R}^d$, we define:
\begin{equation}\label{def-con}\overline{B}_n(x_1,\widetilde{x})\overline{\Phi}^{x_1}(0,\widetilde{x}-x_1\mu)=\Phi^{x_1}(0,\widetilde{x}-x_1\mu)B_n(\widetilde{x}-x_1\mu).\end{equation}
Clearly, $\overline{B}_n(x_1,\widetilde{x})$ are analytic and
$2$-periodic in $\widetilde{x}$. We now
 check that they are  also $2$-periodic in $x_1$.
By $(\ref{almost-1})$ and $(\ref{def-con}),$ we have
\begin{eqnarray*}
&&\overline{B}_n(x_1+2,\widetilde{x})\overline{\Phi}^{x_1}(2,\widetilde{x}-
x_1\mu)\overline{\Phi}^2(0,\widetilde{x}- x_1\mu-2\mu)\\
&=&
\overline{B}_n(x_1+2,\widetilde{x})\overline{\Phi}^{x_1+2}(0,\widetilde{x}-x_1\mu-2\mu)
\\
&=&\Phi^{x_1+2}(0,\widetilde{x}-x_1\mu-2\mu)B_n(\widetilde{x}-x_1\mu-2\mu)\\
&=&\Phi^{x_1}(2,\widetilde{x}- x_1\mu)\Phi^1(0,\widetilde{x}-
x_1\mu-2\mu)B_n(\widetilde{x}-x_1\mu-2\mu)\\&=&\Phi^{x_1}(0,\widetilde{x}-x_1\mu)B_n(\widetilde{x}-x_1\mu)\overline{\Phi}^2(0,\widetilde{x}-
x_1\mu-2\mu)\\
&=&\overline{B}_n(x_1,\widetilde{x})\overline{\Phi}^{x_1}(0,\widetilde{x}-x_1\mu)\overline{\Phi}^2(0,\widetilde{x}-
x_1\mu-2\mu),
\end{eqnarray*}
which means
$\overline{B}_n(x_1+2,\widetilde{x})=\overline{B}_n(x_1,\widetilde{x}),$
and then $\overline{B}_n\in C_{\frac{h_n}{1+|\mu|}}^\omega(2{\T}^d,
SL(2,{\R})).$  By similar reasoning as above, we have

$$\overline{B}_n(\theta_1+t,\widetilde{\theta}+t\mu)\overline{\Phi}^{t}(\theta_1,\widetilde{\theta})=\Phi^{t}(\theta_1,\widetilde{\theta})\overline{B}_n(\theta_1,\widetilde{\theta}).$$
This means $\overline{B}_n$ conjugate $(\omega, A)$ to
$(\ref{al-ref2})$, it concludes that  $(\omega, A)$ is almost
reducible.\end{pf}

\begin{Lemma}
Let $h>0,$ $A\in C^\omega_h(\T^d,sl(2,\mathbb{R}))$, $(\mu,
\mathcal{A})$ is  the corresponding Poincar\'{e} cocycle of
$(\omega,A)$, then $(\mu,\mathcal{A})$ is rotations reducible,  if
and only if $(\omega,A)$ is rotations reducible.
\end{Lemma}

\begin{pf} If $(\omega,A)$ is rotations reducible, the corresponding
Poincar\'{e} cocycle is clearly rotations reducible by definition.

Now we prove the converse part. If $(\mu,\mathcal{A})$ is rotations
reducible, then there exist $h_*<h,$ $B\in
C_{h_*}^\omega(2{\T}^{d-1}, SL(2,\R))$, $\varphi \in
C_{h_*}^\omega({\T}^{d-1}, \R)$ such that
$$\mathcal{A}(\cdot)B(\cdot)=B(\cdot+\mu)R_{\varphi(\theta)}.$$

By Proposition \ref{tech}, there exists  $\rho\in
C_{\frac{h_*}{1+|\mu|}}^\omega({\T}^{d}, \R)$ such that it satisfies
$$ \int_0^1\rho(t,\cdot+t\mu)dt=\varphi(\cdot), $$
which means the quasi-periodic cocycle $(\mu,R_{\varphi(\cdot)})$
can be embedded into the quasi-periodic linear system
$(\omega,\rho(\theta)J).$

We now extend $B(\cdot)$ to the torus $2{\T}^d$ in  the following
way: for $(x_1,\widetilde{x})\in {\R}^d$, we define:
$$\overline{B}(x_1,\widetilde{x})e^{2\pi \int_0^{x_1}\rho(t,\widetilde{x}-x_1 \mu+t\mu)dt J}=\Phi^{x_1}(0,\widetilde{x}-x_1\mu)B_n(\widetilde{x}-x_1\mu).$$

By the same reasoning as in Lemma \ref{almost2},  it can be checked
that
 $$\overline{B}\in C_{\frac{h_*}{1+|\mu|}}^\omega(2{\T}^d,
SL(2,{\R})),$$ and
$$\overline{B}(\theta_1+t,\widetilde{\theta}+t\mu)e^{2\pi
\int_0^{t}\rho(\theta_1+s,\widetilde{\theta}+s\mu)ds
J}=\Phi^{t}(\theta_1,\widetilde{\theta})\overline{B}_n(\theta_1,\widetilde{\theta}),$$
which means $\overline{B}$ conjugate $(\omega, A)$ to
$(\omega,\rho(\theta)J)$, it concludes that  $(\omega, A)$ is
rotations reducible.\end{pf}

\textbf{Proof of Theorem \ref{full-rot}.}

First we recall  results of \cite{HoY}:
\begin{Theorem}\label{hy2}  Let
 $\omega=(1,\alpha)$, $\alpha \in \R\backslash \Q$,
$h>0,$
 $A\in sl(2,\mathbb{R})$ and $F\in C^\omega_{h}(\T^2,sl(2,\R))$. If $$\|F\|_h <C_0 min\{h^{\chi},1\},$$
where $C_0,\chi$ are numerical constants, then  the following
results
 hold:
\begin{enumerate}
\item The system $(\omega,A+F(\theta))$ is almost
reducible.
\item  If the rotation number is Diophantine w.r.t
$\omega$, then it is analytically rotations reducible.
\item  Assume furthermore that
$\beta(\alpha)=0$, then it is analytically reducible.
\end{enumerate}
\end{Theorem}

 In
\cite{HoY}, the transformation converges  on analyticity strips of
width going to zero. As remarked in \cite{HoY}, with minor
modification of the proof, one obtains a strong version of the
almost reducibility below: one can get convergence of the
perturbation on strips of fixed width\footnote{In deed, in order to
prove such strong version of the almost reducibility, one only need
to re-estimate  Lemma $5.2$ of \cite{HoY}.}. If we use such result,
one can easily prove that in Theorem \ref{hy2}, if $2\pi h>\beta>0,$
the rotation number is Diophantine w.r.t $\omega$, then the system
is reducible.

We now finish the proof of Theorem \ref{full-rot}. Suppose
$\|\mathcal {A}-R\|_h=\varepsilon< e^{-\frac{12\pi h}{1+\alpha}},$
then by Theorem \ref{localemb-sl}, there exist $\widetilde{R}\in
sl(2,\R)$ and $\widetilde{F}\in
C^\omega_{h/1+\alpha}(\T^2,sl(2,\R))$ such that the cocycle
$(\alpha,\mathcal {A})$ can be embedded into the quasi-periodic
linear system $(\omega,\widetilde{R}+\widetilde{F})$ with estimate
$\|\widetilde{F}\|_{h/1+\alpha}\leq e^{\frac{6\pi
h}{1+\alpha}}\varepsilon^{\frac{1}{2}}.$ Thus if $$e^{\frac{6\pi
h}{1+\alpha}}\varepsilon^{\frac{1}{2}}< C_0
min\{(\frac{h}{1+\alpha})^{\chi},1\}, $$ which means $$\varepsilon<
\widetilde{C} min\{h^{2\chi},1\}e^{-\frac{12\pi h}{1+\alpha}},$$
then we can apply Theorem \ref{hy2} and Theorem \ref{dis-con} to
finish the proof.\qed

\section{Proof of Theorem \ref{anderson}}

By Aubry duality, we know that the dual operator of $L_{\lambda
V,\alpha, \varphi}$ is $H_{\lambda V,\alpha, \phi}$. The
eigenfunction equation  $H_{\lambda V,\alpha, \phi}x=E x$
corresponds to the Schr\"odinger cocycle $(\alpha,S_E^{\lambda V}).$
We assume that  $$\|\lambda V(\theta)\|_h<\varepsilon< \widetilde{C}
min\{h^{2\chi},1\}e^{-\frac{12\pi h}{1+\alpha}},$$ and  denote by
$$R_{\lambda V,\alpha,\varphi}=\{E|\rot_f(\alpha,S_E^{\lambda
V})=\varphi+\frac{1}{2}\langle k,\alpha\rangle \quad \mbox{mod}
1\},$$
$$\Phi=\{\varphi| \varphi \quad \mbox {is Diophantine w.r.t}\quad \alpha
\},$$ and recall that $\sigma^L_{pp}({\lambda V,\alpha,\varphi})$
and $B_{\lambda V,\alpha,\varphi}$ have been defined in section
\ref{pre-aub}.

The proof of Theorem \ref{anderson} is distinguished into two steps.
First we prove the following:

\begin{Lemma}\label{bloch-redu}
If $\varphi\in \Phi,$ then $R_{\lambda
V,\alpha,\varphi}=\sigma^L_{pp}({\lambda V,\alpha,\varphi}).$
\end{Lemma}

\begin{pf} By Lemma \ref{aubry}, it is sufficient for us to prove
\begin{equation}\label{rot-bloch}
R_{\lambda V,\alpha,\varphi}=B_{\lambda V,\alpha,\varphi} \quad
\mbox{if}\quad \varphi\in \Phi.
\end{equation}

First we prove:
\begin{equation}\label{rot-bloch2}
R_{\lambda V,\alpha,\varphi}\supset B_{\lambda V,\alpha,\varphi}.
\end{equation}
If $E\in B_{\lambda V,\alpha,\varphi},$ then by Lemma $16$ of
\cite{Pui06}, there exists $B\in C^\omega(2\T,SL(2,\R))$ such that
$$B(\theta+\alpha)S_E^{\lambda v}(\theta)B(\theta)^{-1}=\left(
\begin{array}{ccc}
 cos \varphi &  sin \varphi\\
 -sin\varphi &  cos \varphi
 \end{array}\right),$$
thus by Proposition \ref{rot-conj}, we have
$\rot_f(\alpha,S_E^{\lambda V})=\varphi+\frac{1}{2}\langle
k,\alpha\rangle \quad \mbox{mod} 1,$ which proves
$(\ref{rot-bloch2}).$

Then we prove
\begin{equation}\label{rot-bloch1}
R_{\lambda V,\alpha,\varphi}\subset B_{\lambda V,\alpha,\varphi}.
\end{equation}

If $\varphi\in \Phi$ is Diophantine w.r.t $\alpha,$ then
$\rot_f(\alpha,S_E^{\lambda V})$ is also Diophantine w.r.t $\alpha.$
In fact,
$$\|2\varphi-k'\alpha\|_{\R/\Z}\geq \frac{\gamma}{|k'|^\tau},$$
implies that
\begin{eqnarray*}
\|2\rot_f(\alpha,S_E^{\lambda V})-k'\alpha\|_{\R/\Z}&=&
\|2\varphi+k\alpha-k^{'}\alpha\|_{\R/\Z}\\
&\geq& \frac{\gamma}{|k-k^{'}|^\tau} \geq
\frac{(1+|k|)^{-\tau}\gamma}{|k'|^\tau}.
\end{eqnarray*}
Then by Corollary \ref{full-rot}, $(\alpha,S_E^{\lambda V})$ is
rotations reducible, since we assume $2\pi h>(1+\alpha)\beta,$  in
fact we have $(\alpha,S_E^{\lambda V})$ is reducible. So there
exists analytic $B:\T\to SL(2,\R),$ $k\in\Z,$ such that
$$B(\theta+\alpha)S_E^{\lambda V}(\theta)B(\theta)^{-1}=\left(
\begin{array}{ccc}
e^{2\pi i(\varphi+k\alpha)} & 0\\
 0 & e^{-2\pi i(\varphi+k\alpha)}
 \end{array}\right).$$
It follows that  the solution of $H_{\lambda V,\alpha,\phi}x=Ex$ has
quasi-periodic Bloch waves with Floquet exponent $\varphi+k\alpha.$
 Therefore $(\ref{rot-bloch1})$ holds.\end{pf}

Then we prove the following:
\begin{Lemma}
If $\varphi\in \Phi,$ then $\overline{R_{\lambda
V,\alpha,\varphi}}=\sigma^H(\lambda V,\alpha).$
\end{Lemma}

\begin{pf} By $(\ref{rot-bloch})$, we have
$$\overline{R_{\lambda
V,\alpha,\varphi}}=\overline{B_{\lambda V,\alpha,\varphi}}\subset
\sigma^H(\lambda V,\alpha),$$ it is sufficient for us to prove that
$$\sigma^H(\lambda V,\alpha)\subset \overline{R_{\lambda
V,\alpha,\varphi}}.$$

The crucial observation is that when restricted to the spectrum, the
rotation number is strictly monotonic. To simplify the notation, we
let $\rot_f(\alpha,S_E^{\lambda V})=\rot_f(E).$

Since for any $\varphi\in \T,$ the orbit of $\{\varphi+k\alpha
\quad\mbox{mod} 1\}$ is dense in $[0,1]$, we can find $E_n\in
R_{\lambda V,\alpha,\varphi}$ for any $E_0\in \sigma^H(\lambda
V,\alpha),$ such that
$$\rot_f(E_n)\rightarrow \rot_f(E_0),\quad n\rightarrow\infty.$$
Moreover, we can assume $(\rot_f(E_n))_{n\in\Z}$ is monotonic, since
the rotation number is  monotonic, we then have $E_n$ is monotonic
and bounded (the boundness follows from the compactness of
$\sigma^H(\lambda V,\alpha)$). Thus there exists $\widetilde{E}\in
\sigma^H(\lambda V,\alpha),$ such that $E_n\rightarrow
\widetilde{E}.$ By the continuity of the rotational number, we have
$$\rot_f(E_n)\rightarrow \rot_f(\widetilde{E}),\quad
n\rightarrow\infty.$$ Since the rotation number is strictly
monotonic when it restricted to the spectrum $\sigma^H(\lambda
V,\alpha),$ then $\widetilde{E}=E_0,$ it follows that $E_0 \in
\overline{R_{\lambda V,\alpha,\varphi}}.$ \end{pf}

By the last two steps and Lemma \ref{aubry}, we have for any
$\varphi\in \Phi$,
$$\overline{\sigma^L_{pp}({\lambda V,\alpha,\varphi})}=\overline{R_{\lambda
V,\alpha,\varphi}}=\sigma^H(\lambda V,\alpha)=\sigma^L(\lambda
V,\alpha),$$ which means that the long-range operator $L_{\lambda
V,\alpha,\varphi}$ has Anderson Localization. \qed

\textbf{Proof of Theorem \ref{mathieu}:} If we apply Theorem
\ref{anderson} to almost Mathieu operator directly, we can not get
the best estimate, since by the local embedding theorem, we lose
some analyticity, in order to get finer estimate, we need the
following:

\begin{Theorem}\label{hy1} For every
$\alpha \in \R\backslash \Q$, $h>0,$ $\tau>0,$ $\gamma>0,$ let
 $\mathcal {A}\in SL(2,\mathbb{R})$ be such that $\|\mathcal
 {A}-R\|_h<c(h\gamma)^\tau$ for some constant $R$, and
 $\rot_f(\alpha,\mathcal {A})\in DC_\alpha(\tau,\gamma),$ then the system $(\omega,A+F(\theta))$ is analytically
rotations reducible.
\end{Theorem}

\begin{Remark}\label{local-rot}
 The continuous version Theorem  of \ref{hy1} appear in \cite{HoY}.
 In fact,  the proof of this theorem in \cite{HoY} applies essentially unchanged to the
discrete case, thus the same result holds.
\end{Remark}

By Aubry duality, we know that the  almost Mathieu operator is
self-dual, and the dual of $L_{\lambda cos,\alpha, \varphi}$ is
$H_{\lambda cos,\alpha, \phi}$. Suppose
$$\lambda e^{2\pi h}=\|\lambda
cos(2\pi\theta)\|_{h}\leq\varepsilon(h)\ll1$$ is small enough, since
in this case, $\beta(\alpha)>0,$ and it lies in the subcritical
regime, then by Avila's theorem \cite{A}, one can provide (without
any restrictions on the fibered rotation number) a sequence of
conjugacies which put the cocycle arbitrarily close to constants,
(we only lose arbitrary small analyticity strips  of width), so that
Theorem \ref{hy1} eventually can be applied. Thus if the rotation
number is Diophantine w.r.t $\omega$, $2\pi h>\beta(\alpha)>0,$ then
the cocycle is  rotations reducible and consequently reducible: if
$f\in C^\omega_{h}(\T,\R),$ $2\pi h>\beta(\alpha),$ then
$$\varphi(\theta+\alpha)-\varphi(\theta+\alpha)=f(\theta)-\widehat{f}(0)$$
has an analytic solution $\varphi\in
C^\omega_{h-\beta/2\pi}(\T,\R).$ Finally, we can replay the proof of
Theorem \ref{anderson} to finish the proof. \qed

\section*{Acknowledgements}The work was supported by NNSF of China (Grant 10531050),  NNSF
of China (Grant 11031003) and a project funded by the Priority
Academic Program Development of Jiangsu Higher Education
Institutions.  The authors would like  to thank Xuanji Hou, Haken
Eliasson, David Damanik  and Svetlana Jitomirskaya for useful
discussions.

\end{document}